\documentclass[12pt, leqno]{article}

\usepackage{amsmath, amssymb, amscd, verbatim, xspace,amsthm}
\usepackage{latexsym, epsfig, color}

\newcommand{\C}{\ensuremath{{\mathbb{C}}}}
\newcommand{\Z}{\ensuremath{{\mathbb{Z}}}\xspace}

\newcommand{\Q}{\ensuremath{{\mathbb{Q}}}}
\newcommand{\R}{\ensuremath{{\mathbb{R}}}}

\newcommand{\ra}{\rightarrow}

\newcommand\Hom{\operatorname{Hom}}

\newcommand\Aut{\operatorname{Aut}}

\newcommand\im{\operatorname{im}}
\newcommand\Stab{\operatorname{Stab}}
\newcommand\Gal{\operatorname{Gal}}

\newcommand\Ind{\operatorname{Ind}}

\newcommand\f{\mathfrak{f}}
\newcommand\tensor{\otimes}
\newcommand\isom{\cong}
\newcommand\sub{\subset}
\newcommand\tesnor{\otimes}

\newcommand\Disc{\operatorname{Disc}}

\newcommand\Frob{\operatorname{Frob}}

\newcommand{\Nv}{Nv}

\newcommand{\cA}{\mathcal{A}}

\newcommand\bij[2]{\ensuremath{\left\{\parbox{2.5 in}{#1}\right\} \longleftrightarrow \left\{\parbox{2.5 in}{#2}\right\}}}

\newtheorem{proposition}{Proposition}[section]
\newtheorem{theorem}[proposition]{Theorem}
\newtheorem{corollary}[proposition]{Corollary}

\newtheorem{lemma}[proposition]{Lemma}

\theoremstyle{remark}
\newtheorem{remark}[proposition]{Remark}

\newcommand{\E}{\mathcal{E}(C)}

\newcommand\Div{\operatorname{Div}}
\newcommand\Sp{\operatorname{Sp}}
\renewcommand\Re{\operatorname{Re}}

\usepackage{fullpage}

\newenvironment{notation}{\vspace{2 ex}{\noindent{\bf Notation. }}}{\vspace{2 ex}}
\renewcommand\bij[2]{\ensuremath{\left\{\parbox{2 in}{#1}\right\} \longleftrightarrow \left\{\parbox{2 in}{#2}\right\}}}
\newenvironment{proof2}{\noindent {\it Proof }}{$\Box$ \vspace{2 ex}}

\title{On the probabilities of local behaviors in abelian field extensions}

\author{Melanie Matchett Wood}

\begin{document}

\maketitle

\begin{abstract}
For a number field $K$ and a finite abelian group $G$,
we determine the probabilities of various local completions of a random
$G$-extension of $K$ when extensions are ordered by conductor.
In particular, for a fixed prime $\wp$ of $K$, we determine
the probability that $\wp$ splits into $r$ primes in a random
$G$-extension of $K$ that is unramified at $\wp$.  
We find that these probabilities
are nicely behaved and mostly independent. 
This is in analogy to Chebotarev's density theorem, which
gives the probability that in a fixed extension a random prime of $K$
splits into $r$ primes in the extension.   
We also give the asymptotics for the number of $G$-extensions with bounded conductor.
In fact, we give a class of extension invariants, including
conductor, for which we obtain the same counting and probabilistic results.  In contrast, we prove that
that neither the analogy with the Chebotarev probabilities nor the independence of probabilities holds
when extensions are ordered by discriminant. 
\end{abstract}

\section{Introduction}\label{intro}
 
Given a finite Galois extension $L/\Q$ with
Galois group $G$, and a rational prime $p$, what is the probability that
$p$ splits completely in $L$?  If we fix $L$ and vary 
$p$, the Chebotarev density theorem tells us what proportion of
primes have any given splitting behavior.  However, we can alternatively fix $p$
(and $G$), and study the probability that $p$ splits a certain way in a random $L$
with $\Gal(L/\Q)\isom G$.  We ask whether the probabilities of the unramified splitting types
 are in the proportions we expect from the Chebotarev density theorem.
We also ask if the probabilities are independent at different primes $p$. In fact, we shall ask more refined questions and study the probabilities of various local
 $\Q_v$-algebras $L_v:=L\tensor_\Q \Q_v$ at a place $v$ of $\Q$.  
These questions have recently been asked by Bhargava \cite[Section 8.2]{B:MF} and have come up naturally
in the work counting extensions of $\Q$ with a given Galois group (see \cite{Cohen1}, \cite{Cohen2}, \cite{Taylor}, \cite{Wright}, and Section~\ref{SS:hist}).
In this paper, we
answer these refined questions for abelian $G$. For the rest of this paper, we fix a finite abelian group $G$.

We define a \emph{$G$-extension of} a field $K$ to be a Galois extension $L/K$ with 
an isomorphism $\phi : G \ra \Gal(L/K)$.  
An isomorphism of two $G$-extensions $L$ and $L'$ is given by 
an isomorphism $L\ra L'$ of $K$-algebras that respects the $G$-action on $L$ and $L'$.
Let $E_G(K)$ be the set of isomorphism classes of $G$-extensions of $K$.
Given a finite set $S$ of places of $\Q$, and  a $\Q_v$-algebra $T_v$ for each $v\in S$, we use $T$ to denote the collection of
all the choices $T_v$.  We define the \emph{probability of $T$} as follows: 
\begin{equation}\label{E:probdef}
 \Pr(T)=\lim_{X\ra\infty}\frac{\#\{\textrm{ $L\in E_G(\Q)$ $\mid$ $L_v\isom T_v$ for all $v\in S$ and $\f(L)<X$}\}}
{\#\{\textrm{$L\in E_G(\Q)$  $\mid$ $\f(L)<X$}\}},
\end{equation}
where $\f(L)$ is the finite conductor of $L$ over $\Q$.  We can analogously define the probability
of one local algebra $T_v$, or of a splitting type of a prime.

Given a $G$-extension $L$ of $\Q$, every $L_v$ is of the form $M^{\oplus r}$, where
$M$ is a field extension of $\Q_v$ with Galois group $H$, and $H$ is a subgroup of $G$ of index $r$.
The first twist in this story is that some $M^{\oplus r}$ of this form \emph{never}
occur as $L_v$.  For example, when $G=\Z/8\Z$, it is never the case that $L_2/\Q_2$ is unramified of degree 8.
This means we cannot expect unramified splitting types to occur in the proportions suggested by the Chebotarev density theorem.
Wang \cite{Wang}, in a correction to work of Grunwald \cite{Grunwald}, completely determined which 
 local algebras occur.
The only obstruction is that for even $|G|$, some $\Q_2$ algebras do not occur
as $L_2$ for any $G$-extension $L$.
Call these $\emph{inviable}$ $\Q_2$-algebras 
(and all other $M^{\oplus r}$ of the above form \emph{viable})
and note the characterization implicitly depends on $G$.
Once one knows which local algebras can occur, it is natural to ask how often they occur.
We answer that question in the following theorem.
\begin{theorem}\label{T:pre}
 Let $v$ be a place of $\Q$. Let
$M$ and $M'$ be field extensions of $\Q_v$ with Galois groups $H$ and $H'$ that are subgroups of $G$ of index $r$ and $r'$, 
respectively.
Then, unless $v=2$ and at least one of $M^{\oplus r}$, $M'^{\oplus r'}$ is inviable,
$$
\frac{\Pr(M^{\oplus r})}{\Pr(M'^{\oplus r'})}=\frac{\frac{|\Hom_0(H,G)|}{\f(M)}}{\frac{|\Hom_0(H',G)|}{\f(M')}},
$$
where $\Hom_0(E,G)$ denotes the set of injective homomorphisms from $E$ to $G$.  
The conductor $\f(M)$ is viewed as an element of $\Q$.
\end{theorem}

We will refer to the density of primes with a given splitting type in a fixed $G$-extension 
as the \emph{Chebotarev probability} of that splitting type.  We compare Theorem~\ref{T:pre} to the
Chebotarev density theorem in the following corollary.

\begin{corollary}\label{C:QCheb}
The probability of a fixed rational prime $p$ (not $2$ if $|G|$ is even) splitting 
into $r$ primes in a random $L\in E_G(\Q)$, given that
$p$ is unramified, is the same as the Chebotarev probability
of a random rational prime $p$ splitting into $r$ primes in
 a fixed $L\in E_G(\Q)$.
\end{corollary}

In fact, it follows from Theorem~\ref{T:pre} that when $|G|$ is even and $p=2$ the probabilities of 
viable splitting types in a random $G$-extension occur in the same proportions as they occur in the
Chebotarev density theorem for a fixed extension and random prime.
Of course, one contrast to the Chebotarev probabilities is that for a fixed $p$ and a random $G$-extension $L$, the prime
$p$ will be ramified with positive probability.  In
 this paper, we also determine the independence of the local probabilities computed in Theorem~\ref{T:pre}, leading to the
following result.

\begin{theorem}\label{T:Qind}
For any finite set $S$ of places of $\Q$ and any choice of local $\Q_v$-algebras $T_v$ for $v\in S$, 
the events $T_v$ are independent.  
\end{theorem}

One may ask whether we obtain the same result if we count the $G$-extensions in other ways,
for example by replacing the conductor by the discriminant, by an Artin conductor, or by the product of the ramified primes.
In fact, in Section 2 we prove a stronger version of Theorem~\ref{T:pre} which replaces the conductor
with any function satisfying a certain \emph{fairness} hypothesis (defined in Section 2),
 which is satisfied by the conductor, some Artin conductors, and the product of the ramified primes.  In Section~\ref{S:AC} we give examples of some Artin conductors that are fair.
The discriminant is fair only when $G$ has prime exponent.
Much work has been done to study the asymptotics of the number of extensions with bounded discriminant
and having Galois closure with a specified Galois group (see \cite{Cohen1} and \cite{Cohensur} for surveys).
These asymptotics were determined completely for abelian Galois groups by M\"{a}ki \cite{MakiDisc}.
M\"{a}ki \cite{MakiCond} also has determined
 the asymptotics of the number of extensions with fixed abelian Galois group and bounded conductor.
  In Section~\ref{S:constant}, we give the 
asymptotics of the number of $G$-extensions with bounded conductor 
(or any fair counting function)
for a finite abelian group $G$.  Our result is a generalization of M\"{a}ki's work \cite{MakiCond},
in that we can replace the conductor by other fair counting functions and that we give the result over an arbitrary base number field
 (see Section~\ref{SS:other}).  We also give the constant in the asymptotic more explicitly than it appears in \cite{MakiCond}.

For degree $n$ extensions having Galois closure with Galois group $S_n$, it is known that  when counting
by discriminant for $n=2,$ $3,$ $4$, and $5$, local completions $L_v$ show up with probability proportional to
$\frac{1}{|\Aut_{\Q_v} L_v|}\frac{1}{|\Disc L_v|}$ (see \cite{S3}, \cite{S4}, and \cite{S5} for the computation of the probabilities,
and \cite{B:MF} for this interpretation). 
We will see after Corollary~\ref{C:probs} how to interpret our probabilities in closer analogy to the results in \cite{S3}, \cite{S4}, and \cite{S5}.
However, it turns out that counting abelian extensions by discriminant does not lead to such nice local probabilities.
This was observed by Wright \cite{Wright}, in his work on counting abelian extensions asymptotically by discriminant.
Let the \emph{discriminant probability} be defined as in Equation~\eqref{E:probdef} but with the conductor replaced by the 
absolute value of the discriminant.  We call two events \emph{discriminant independent} if they are
independent with the discriminant probability.
Wright showed that all viable $\Q_v$-algebras occur with positive discriminant probability,
and noted that when $G$ has prime exponent, the relative probabilities of local extensions are simple expressions.
(Wright actually works over an arbitrary global field with characteristic not dividing $|G|$; 
in Section~\ref{SS:other} of this paper we describe our work over an arbitrary number field.)
When $G=\Z/4\Z$, Wright notes that the ratio of the discriminant probability of $\Q_p^{\oplus 4}$ to the
the discriminant probability of the unramified extension of $\Q_p$ of degree 4 is an apparently very complicated expression.
In Section~\ref{S:disc}, we prove the following propositions in order to show that the discriminant probability
analogs of Corollary~\ref{C:QCheb} or Theorem~\ref{T:Qind} do not hold.

\begin{proposition}
Let $p$, $q_1,$ and $q_2$ be primes with $q_i\equiv 1\pmod{p^2}$ for $i=1,2$.
 Then $q_1$ ramifying and $q_2$ ramifying in a random $\Z/p^2\Z$-extension
are not discriminant independent. 
\end{proposition}

The Chebotarev probability that a 
random prime splits completely in a fixed $\Z/9\Z$-extension is $\frac{1}{9}$.
However, we have the following.

\begin{proposition}\label{P:disccheb}
Let $q=2$, $3,$ $5,$ $7,$ $11,$ or $13$.
Given that $q$ is unramified, the discriminant probability that $q$ splits completely
in a random $\Z/9\Z$-extension is strictly less than $\frac{1}{9}$.
\end{proposition}

For comparison, in the above two cases we have that the (conductor) 
probabilities are independent, and the (conductor) probability is $\frac{1}{9}$, respectively.

\subsection{Other base fields}\label{SS:other}
Of course, we can ask all of the same questions when $\Q$ is replaced by an arbitrary number field $K$, and we now fix a number field $K$.
However, for arbitrary number fields there is a further twist in this story.
Given $G$, it is possible that the $K_v$-algebra $T_v$ and the $K_{v'}$-algebra $T'_{v'}$ both occur
from global $G$-extensions, but never occur simultaneously (see \cite{GWex}).
This suggests that we should not expect $T_v$ and $T'_{v'}$ to be independent events.
However, given obstructions of this sort, which were completely determined in \cite{Wang} (or see \cite[Chapter 10]{AT}), we have the best possible
behavior of the local probabilities.  We shall need more precise language to clearly
explain this behavior.

The local $K_v$-algebras coming from $L$ have  structure that we have  so far ignored; namely, they have
a $G$-action coming from the global $G$-action.
Given a field $F$, a \emph{$G$-structured $F$-algebra} is an \'{e}tale $F$-algebra $L$ of degree $|G|$  with an
inclusion $G\hookrightarrow \Aut_F(L)$ of $G$ into the $F$-algebra automorphisms of $L$, such that $G$ acts transitively on the 
idempotents of $L$. 
An isomorphism of two $G$-structured  $F$-algebras $L$ and $L'$  is an
$F$-algebra isomorphism $L{\ra} L'$ such that the induced map
$\Aut_F(L) \ra \Aut_F(L')$ restricts to the identity on $G$.  
If we have a $G$-extension $L$ of $K$, for each place $v$ of $K$, then we have a $G$-structured
$K_v$-algebra $L_v= L \tesnor_K K_v$, where $G$ acts on the left factor.  
Given a subgroup $H$ of $G$, and an $H$-extension $M$ of $K_v$, we can 
form the induced $G$-structured $K_v$-algebra $\Ind_H^G M$ via the usual construction of an induced representation,
which will have a natural structure of an \'{e}tale $K_v$-algebra.  All $G$-structured
$K_v$-algebras coming from $G$-extensions $L$ of $K$ are of the form $\Ind_H^G M$.
So we can ask an even more refined question, at all places,
 about the probability of a certain $G$-structured $K_v$-algebra.
We let $\f(L)$ be the norm from $K$ to $\Q$ of the conductor of $L/K$ (or of the conductor of $L/K_v$, viewed as an ideal of $K$).
Let $S$ be a finite set of places of $K$, and let $\Sigma$ denote
a choice $\Sigma_v$ of $G$-structured $K_v$-algebra for each $v\in S$,
which we refer to as a (\emph{local}) \emph{specification}.
We can then define probabilities as in Equation~\eqref{E:probdef}, replacing $E_G(\Q)$ with $E_G(K)$.

If there exists a $G$-extension $L/K$ such that $L_v\isom \Sigma_v$ for all $v\in S$, then we call $\Sigma$ \emph{viable}
and otherwise we call it \emph{inviable}.
The question of which specifications are viable has been completely answered (see \cite[Chapter 10]{AT}).
There is a set $S_0$ of places of $K$ (depending on $G$, all dividing 2, and empty if $|G|$ is odd) and a finite list 
$\Sigma(1),\dots,\Sigma(\ell)$ of local specifications on $S_0$ such that a local specification $\Sigma$ on $S$ is viable if and only if 
either $S_0\not\sub S$ or $\Sigma$ restricts to some $\Sigma(i)$ on $S_0$.  (We give $S_0$ explicitly in Section~\ref{S:anal}.)
In other words, whether a specification on $S$ is viable depends only on its specifications at places in $S_0$,
and if a specification does not include specifications at all places in $S_0$ then it is viable. 

Now we will build a model for the expected probabilities of local specifications.
Let $\Omega=\prod_{v\textrm{ place of } K} \{ \textrm{isom. classes of $G$-structured $K_v$-algebras} \}$.  
For a local specification $\Sigma$, let $\tilde{\Sigma}=\{x\in \Omega \mid  x_v \isom \Sigma_v \textrm{ for all } v\in S\}$. Let $A=\bigcup_{i=1}^{\ell}\tilde{\Sigma}(i)$, where $\Sigma(i)$ are as in the above paragraph
in the condition for a local specification to be viable.  So for a specification $\Sigma$
on $S$, we have that $\tilde{\Sigma}\cap A$ is non-empty if and only if $\Sigma$ is viable, and in fact
$$\tilde{\Sigma}\cap A=\{ \tilde{\Sigma'} \textrm{ local specification on } S\cup S_0 \mid \Sigma' \textrm{ viable and restricts to } \Sigma \textrm{ on $S$}\}.$$
The $\tilde{\Sigma}_v$ generate an algebra of subsets of $\Omega$.  We can define a finitely additive probability
 measure $P$ on
this algebra by specifying that 
\begin{enumerate}
 \item $\displaystyle{\frac{P(\tilde{\Sigma}_v)}{P(\tilde{\Sigma'}_v)}=\frac{\frac{1}{\f(\Sigma_v)}}{\frac{1}{\f(\Sigma'_v)}}}$
for all  $G$-structured $K_v$-algebras $\Sigma_v$ and $\Sigma_v'$ 
\item $\tilde{\Sigma}_{v_1},\dots,\tilde{\Sigma}_{v_s}$ at pairwise distinct places $v_1,\dots,v_s$, respectively, are independent.
\end{enumerate}
We might at first hope
that $P$ is a model for the probabilities of local specifications in the space of $G$-extensions.  However,
once we know that some specifications never occur, including combinations of occurring specifications, the best we can hope for is the following, which we prove in Section~\ref{S:anal}.
\begin{theorem}\label{T:best}
 For a local specification $\Sigma$ on a finite set of places $S$,
$$
\Pr(\Sigma)=P(\tilde{\Sigma} | A).
$$
\end{theorem}
\begin{corollary}\label{C:probs}
 If $S$ is a finite set of places of $K$ either containing $S_0$ or disjoint from $S_0$, and
$\Sigma$ and $\Sigma'$ are viable local specifications on $S$ then
$$
 \frac{\Pr(\Sigma)}{\Pr(\Sigma')}=\frac{\prod_{v\in S} \frac{1 }{\f(\Sigma_v)}
  }{\prod_{v\in S} \frac{1}{\f(\Sigma'_v)}}.
$$
\end{corollary}

All $G$-structured algebras have $|G|$ automorphisms
(Proposition~\ref{P:autG}), and so for $v$ not in $S_0$, we can also say that the probability of $\Sigma_v$ is proportional to
$\frac{1 }{|\Aut(\Sigma_v)|\f(\Sigma_v)}$.

\begin{corollary}\label{C:Cheb}
The probability of a fixed prime $\wp$ of $K$ (not in $S_0$) splitting into $r$ primes in a random $L\in E_G(K)$, given that
$\wp$ is unramified, is the same as the Chebotarev probability
of a random prime $\wp$ of $K$ splitting into $r$ primes in a fixed $L\in E_G(K)$.
\end{corollary}
\begin{corollary}\label{C:ind}
 If $S_1, \dots, S_t$ are pairwise disjoint finite sets of places of $K$, and each $S_i$
either contains $S_0$ or is disjoint from $S_0$, then local specifications
$\Sigma^{(i)}$ on $S_i$ are independent.
\end{corollary}

Theorem~\ref{T:best} says that the probabilities of local specifications of random $G$-extensions are exactly as in a model
with simple and independent local probabilities, but restricted to a subspace corresponding to the viable specifications on $S_0$.  As when $K=\Q$, we prove Theorem~\ref{T:best}
and its corollaries as a special case of analogous results (see Theorem~\ref{T:prob}) for more general ways of
counting extensions than by conductor.

\subsection{History of the problem and previous work}\label{SS:hist}
The results mentioned above 
of Davenport and Heilbronn (\cite{S3}) and Bhargava (\cite{S4} and \cite{S5})
are a major motivation of this work.
These results show that the local behaviors of
random degree $n$ extensions of $\Q$ whose Galois closure has Galois group $S_n$
have nice discriminant probabilities and are discriminant independent,
when $n=3$, $4,$ or $5$.
The work of Datskovsky and Wright \cite{DW} generalizes that of Davenport and Heilbronn
(the case $n=3$) to
an arbitrary base field.

Taylor \cite{Taylor} proves the result of our Corollary~\ref{C:Cheb} in the special case 
that $G=\Z/n\Z$, and assuming that if $2^g\mid n$ then $K$ contains the $2^g$th roots of unity
(in which case $S_0$ is empty).  Taylor attributes the question of the distribution of splitting types
of a given prime in random $G$-extensions to Fr\"{o}hlich, who was motivated by the work of Davenport and Heilbronn \cite{S3}.
Wright \cite{Wright} proves an analog of Corollary~\ref{C:probs} 
for discriminant probability in the case that $G=(\Z/p\Z)^b$ for $p$ prime
and $|S|=1$, and for these $G$ the discriminant is a fixed power of the conductor, and thus
discriminant probability is the same as conductor probability.
Wright \cite{Wright} suggests that his methods for counting abelian extensions by discriminant could be combined with the methods of Taylor to count abelian extensions by conductor.  In this paper, we follow this suggestion
and incorporate methods of both Wright and Taylor
along with some new ideas.  We implicitly count abelian extensions by conductor (and give this result in Section~\ref{S:constant}),
but are focused on the probabilities of local behaviors.

In the work of counting extensions
whose Galois closure has some fixed Galois group, it has been often suggested that it is natural to also
count such extensions with fixed local behavior 
(for example, in the work of Cohen, Diaz y Diaz, and Olivier \cite{Cohensur} for the group $D_4$, 
the heuristics of Malle \cite[Remark 1.2]{Malle} for general groups, and
in the general surveys \cite{Cohen1} and \cite{Cohen2}).  Some authors have also considered these questions
when one replaces field extensions with polynomials, and counts with a natural density on the polynomials
(see \cite{DC1}, \cite{DC2}, \cite{DC3}, and \cite{vdP}).  

Theorems~\ref{T:pre}, \ref{T:Qind}, and \ref{T:best} and their corollaries are all new (except in
the special cases mentioned above), but the proofs use many techniques that come from the work of 
Taylor and Wright.  
Some new techniques are required to compute the probabilities exactly
in the case of non-cyclic $G$ and for more general ways of counting extensions.
An important new ingredient is the consideration of the probabilities
of $G$-structured $K_v$-algebras (and not just $K_v$-algebras), which not only
allows us to give more refined probabilities but allows us to state Theorem~\ref{T:best}.
One of the central contributions of this paper is the formulation of 
Theorem~\ref{T:best}, which makes precise the idea that the probabilities are as well-behaved as possible
in light of the non-occurrence of certain local extensions (see \cite{Wang} and \cite{AT}).
For abelian groups $G$, we study for the first time the
probabilities when more than one local behavior is specified and
the independence of these local probabilities.
Our results are for all base number fields $K$, all finite abelian groups $G$, and for many ways of counting extensions
(see the definition of \emph{fair} in Section~\ref{S:anal}) including by conductor.

\subsection{Outline of the paper}
In Section~\ref{S:anal}, we define counting functions and fairness, and prove our main theorems.  
The proof of our main theorems involves making a Dirichlet series generating function for the extensions we are counting,
relating it to $L$-functions whose analytic behavior is known, using standard Tauberian theorems to deduce
asymptotic counting results, and using fairness to express the desired probabilities in a simple form.
In Section~\ref{S:constant}, we give the asymptotic number of $G$-extensions with a given invariant (such as conductor) bounded.
We give an explicit Euler product for the constant in this asymptotic result.
In Section~\ref{S:disc}, we prove that when counting by discriminant, the local probabilities do
not have the same nice behavior as in the conductor case.  In Section~\ref{S:AC}, we give some examples of fair Artin conductors.
In Section~\ref{S:further}, we discuss the further questions that this work motivates.

\section{Statement and proof of the main theorem}\label{S:anal}

 In this section, we prove a generalization of Theorem~\ref{T:pre} and Theorem~\ref{T:best} for more general ways of counting 
$G$-extensions
than by conductor. 
First, in Subsection~\ref{SS:count}, we will define the acceptable ways of counting $G$-extensions.
Then, in Subsection~\ref{SS:state}, we state Theorem~\ref{T:prob}
(our generalization of Theorem~\ref{T:pre} and Theorem~\ref{T:best}) and deduce several corollaries.
In Subsection~\ref{SS:G}, we relate $G$-structured algebras to Galois representations.  In Subsection~\ref{SS:Euler},
we define a generating function counting $G$-extensions satisfying a local specification $\Sigma$ and express
this generating function as a sum of Euler products.  In Subsection~\ref{SS:poles}, we 
state three lemmas about the analytic behavior of these Euler products,
and then use the standard Tauberian analysis to determine the asymptotic behavior of the coefficient sums
of the generating function from the rightmost poles.  From this asymptotic behavior we deduce Theorem~\ref{T:prob}.
In Subsection~\ref{S:cont}, we prove the three lemmas stated in Subsection~\ref{SS:poles}.
The method in Subsection~\ref{SS:Euler} is very similar to that of Wright \cite{Wright} and some of 
the methods in Subsection~\ref{S:cont} are motivated by those of Taylor \cite{Taylor}.

\subsection{Counting functions and fairness}\label{SS:count}
 We fix a finite abelian group $G$ and a number field $K$.  Let $n=|G|$.  Let $c_G : G \ra \Z_{\geq 0}$ be a function such that
1) $c_G(g)=0$ if and only if $g=1$ and 2)  if $e$ is relatively prime to the order of $g\in G$, then $c_G(g^e)=c_G(g)$.  For all places $v$ dividing $n$ or infinite,
let
$c_v : \{\mbox{isom. classes of $G$-structured $K_v$-algebras}\} \ra \Z_{\geq 0}$ be an arbitrary function.
From these functions $c_G$ and the $c_v$, we define
$c : \{\mbox{isom. classes of $G$-structured $K_v$-algebras}\} \ra \Z_{\geq 0}$
by 
$$
c(\Sigma_v)=
\begin{cases}
 {c_G(y_v)} &\text{if } v\nmid n\infty, \text{ where }\Sigma_v=\Ind^G_H M, \text{ and $M/K_v$ a field extension,} \\
&  \text{            and } y_v \text{ is any generator of tame inertia in }\Gal(M/K_v) \subset G;\\
{c_v(\Sigma_v)} &\text{if } v\mid n\infty.
\end{cases}
$$
We then define an invariant
$C$ of $G$-extensions
by the product $C(L)=\prod_{v} \Nv^{c(L_v)}$ over places of $K$, where $\Nv$ is $N_{K/\Q} v$ at finite places and by convention 1 at infinite places. 
 We call such a $C$, determined by components $c_G$ and the $c_v$, a \emph{counting function}. Let $m=\min_{g\in G\setminus \{1\}} c_G(g)$ and let
$\mathfrak{M}=c_G^{-1}(m)$.  
Let $G_r=\{x\in G \mid x^r=1\}$.

A counting function is \emph{fair} if for all $r$, we have that
$\mathfrak{M}\cap G_r$ generates $G_r$.
The norms to $\Q$ of the conductor and  of the product of ramified primes of an extension are both fair counting functions
with $m=1$ and $\mathfrak{M}=G\setminus \{1\}$.  The discriminant is a counting function, but it is not fair unless $G$
has prime exponent.  For example, when $G=\Z/p^2\Z$, for the discriminant we have $\mathfrak{M}=p\Z/p^2\Z$.
In Section~\ref{S:AC}, we give some examples of fair Artin conductors.

\subsection{Statement of the main theorem and corollaries}\label{SS:state}
We define the $C$-probability, $\Pr_C$, by replacing $\f$ with $C$ in Equation~\eqref{E:probdef}.
(Note that $C(L)<X$ implies that $L$ is unramified at all primes
larger than $nX$, and so there are only finitely many such extensions.)
As in the definition of $P$ in the introduction, we define $P_C$
on the algebra of subsets of 
 $\Omega=\prod_{v\textrm{ place of } K} \{ \textrm{isom. classes of $G$-structured $K_v$-algebras} \}$
generated by the $\tilde{\Sigma}_v$ by specifying
\begin{enumerate}
 \item $\displaystyle{\frac{P_C(\tilde{\Sigma}_v)}{P_C(\tilde{\Sigma'}_v)}=\frac{\Nv^{-c(\Sigma_v)/m}}{\Nv^{-c(\Sigma'_v)/m}}}$
for all $G$-structured algebras $\Sigma_v$ and $\Sigma_v'$ and
\item $\tilde{\Sigma}_{v_1},\dots,\tilde{\Sigma}_{v_s}$ at pairwise distinct places $v_1,\dots,v_s$ are $C$-independent.
\end{enumerate}

Let $\eta_i=\zeta_{2^i}+\zeta_{2^i}^{-1}$, where $\zeta_{2^i}$ is a primitive $2^i$th root of unity.
Let $s$ be maximal such that $\eta_s\in K$.  If $2^{s+1}$ does not divide the exponent of $G$, 
then let $S_0=\emptyset$.  Otherwise, let $S_0$
be the set of primes $\wp$ of $K$ dividing 2 such that none of $-1$, $2+\eta_s$ and $-2-\eta_s$
are squares in $K_\wp$.
Recall that there is a list
$\Sigma(1),\dots,\Sigma(\ell)$ of local specifications on $S_0$ such that a local specification $\Sigma$ on $S$ is viable if and only if 
either $S_0\not\sub S$ or $\Sigma$ restricts to some $\Sigma(i)$ on $S_0$ (see \cite[Chapter 10]{AT}). We have defined $A=\bigcup_{i=1}^\ell \tilde{\Sigma}(i)$.  If $S_0$ is empty, then all local specifications are viable and $A$ is the total space $\Omega$.
In this section, we prove the following theorem, of which Theorems~\ref{T:pre} and \ref{T:best} are special cases.
\begin{theorem}\label{T:prob}
 For a local specification $\Sigma$ on a finite set of places $S$ and a fair counting function $C$,
$$
{\Pr}_C(\Sigma)=P_C(\tilde{\Sigma} | A).
$$
\end{theorem}

Now, we will prove several corollaries of Theorem~\ref{T:prob}.  Corollaries~\ref{C:QCheb}, \ref{C:probs}, \ref{C:Cheb}, and
\ref{C:ind} from the introduction are just the following corollaries when $C$ is the norm to $\Q$ of the conductor.
Theorem~\ref{T:Qind} follows from Corollary~\ref{C:ind}.
\begin{corollary}\label{C:probs2}
 If $S$ is a finite set of places of $K$ either containing $S_0$ or disjoint from $S_0$, and
$\Sigma$ and $\Sigma'$ are viable local specifications on $S$ then
$$
 \frac{{\Pr}_C(\Sigma)}{{\Pr}_C(\Sigma')}=\prod_{v\in S} \frac{\Nv^{-c(\Sigma_v)/m}}{\Nv^{-c(\Sigma'_v)/m}}.
$$
\end{corollary}
\begin{proof}
 If $S$ is disjoint from $S_0$, then since $A$ only includes specifications on $S_0$, we have that
$\tilde{\Sigma}$ and $\tilde{\Sigma}'$ are each $P_C$-independent from $A$ in $\Omega$.
Thus ${\Pr}_C(\Sigma)=P_C(\tilde{\Sigma}|A)=P_C(\tilde{\Sigma})$, and similarly for $\Sigma'$.
If $S\supset S_0$, then since $\Sigma$ is viable, $\tilde{\Sigma}\subset A$.  Thus,
${\Pr}_C(\Sigma)=P_C(\tilde{\Sigma}|A)=P_C(\tilde{\Sigma})/P_C(A)$, and similarly for $\Sigma'$.
\end{proof}
\begin{corollary}\label{C:Cheb2}
The $C$-probability of a fixed prime $\wp$ of $K$ (not in $S_0$) splitting into $r$ primes in a random $L\in E_G(K)$, given that
$\wp$ is unramified, is the same as the Chebotarev probability
of a random prime $\wp$ of $K$ splitting into $r$ primes in a fixed $L\in E_G(K)$.
\end{corollary}
\begin{proof}
 The number of $\Sigma_\wp$ that give $\wp$ unramified and splitting into $r$
primes is the number of order $|G|/r$ elements of $|G|$.  (This can be seen, for example, from Lemma~\ref{L:corr}.)
Thus 
$$
\frac{{\Pr}_C(\wp\textrm{ splits unramified into $r$ primes})}{{\Pr}_C(\wp\textrm{ splits unramified into $r'$ primes})}=
\frac{\textrm{number of order $|G|/r$ elements of $|G|$}}{\textrm{number of order $|G|/r'$ elements of $|G|$}},
$$
which agrees with the Chebotarev probabilities.
\end{proof}

\begin{corollary}\label{C:ind2}
 Let $S_1, \dots, S_t$ be pairwise disjoint finite sets of places of $K$, and suppose each $S_i$
either contains $S_0$ or is disjoint from $S_0$. 
(For example, if $|S_0|$ is 0 or 1, then this is always the case.)
 Then local specifications
$\Sigma^{(i)}$ on $S_i$ are $C$-independent. 
\end{corollary}
\begin{proof}
 If $S_0$ is empty, then $A=\Omega$, and this corollary is clear.  Otherwise, first suppose some $S_i$,
say $S_1$, contains $S_0$.
If $\Sigma^{(1)}$ is inviable, then ${\Pr}_C(\Sigma^{(1)})=0$ and otherwise we have
${\Pr}_C(\Sigma^{(1)})=P_C(\tilde{\Sigma}^{(1)})/P_C(A)$.  For
  $i\ne 1$ we have ${\Pr}_C(\Sigma^{(i)})=P_C(\tilde{\Sigma}^{(i)})$,
 as in the proof of Corollary~\ref{C:probs2}.
Let $\Sigma$ be the local specification that is union of the $\Sigma^{(i)}$.
If 
$\Sigma^{(1)}$ is inviable then ${\Pr}_C(\Sigma)=0$, and otherwise
 ${\Pr}_C(\Sigma)=P_C(\tilde{\Sigma})/P_C(A)=\frac{\prod_i P_C(\tilde{\Sigma}^{(i)})}{P_C(A)}=\prod_i \Pr_C({\Sigma}^{(i)})$.

If, on the other hand, no $S_i$ contains $S_0$, then we have
${\Pr}_C(\Sigma^{(i)})=P_C(\tilde{\Sigma}^{(i)})$ for all $i$ 
and $\tilde{\Sigma}$ is $C$-independent from $A$.  Thus ${\Pr}_C(\Sigma)=P_C(\tilde{\Sigma})=\prod_i \Pr_C(\tilde{\Sigma}^{(i)})=\prod_i \Pr_C({\Sigma}^{(i)})$.
\end{proof}

\begin{notation}
We let $n=|G|$ and write $G\isom \Z/n_1 \times \dots \times \Z/n_k$.  For the rest of Section~\ref{S:anal} we use additive notation for $G$.
For all positive integers $m$, we choose compatible primitive $m$th roots of unity $\zeta_m$ such that
if $m' | m$, then $\zeta_{m'}=\zeta_m^{m/m'}$.
  Let $J$ be the group of id\`{e}les of $K$.  For a map
$\chi$ from $J$, we denote by $\chi_v$ the restriction of $\chi$ to $K_v^\times$.
Let $o_v$ be the ring of integers of $K_v$.
Let $J_S$ be the group of id\`{e}les which have components in $o_v^\times$ for all places $v\not \in S$.
In this paper, when we write a map from the id\`{e}les, id\`{e}le class group, or $K_v^\times$ to a finite group (e.g. $\chi : J \ra G$), it will
always mean a continuous homomorphism (for the discrete topology on the range).
\end{notation}
 
\subsection{$G$-structured algebras and Galois representations}\label{SS:G}
Recall that $G$ is a finite abelian group.  The following two results are fairly standard, but we include them here for completeness.
\begin{lemma}\label{L:corr}
For a field $F$, there is a one to one-correspondence 
\begin{align*}
 \bij{\textrm{isomorphism classes of $G$-structured} F\textrm{-algebras}}{\textrm{continuous
homomorphisms\\ $G_F \ra G$}},
\end{align*}
where $G_F$ is the Galois group of a separable closure of $F$ over $F$.
In this correspondence, $G$-extensions correspond to surjective homomorphisms.
\end{lemma}

\begin{proof}
Given a $G$-structured  $K$-algebra $L$ with $G\subset\Aut_K(L)$,
we consider the stabilizer $\Stab\subset G$ of one of the fields $L_0$ that is a direct summand of $L$.
  We have a morphism
$\Stab\ra \Gal(L_0/K)$.  Since $G$ is transitive on the idempotents of $L$ and abelian,
this is an injection.  Since $G$ is transitive on the idempotents,
we see that all the fields that are direct summands of $L$ are isomorphic, and thus 
$|\Stab|=[L_0:K]$.  Thus we have $\Stab\ra \Gal(L_0/K)$ is an isomorphism.
Its inverse gives $G_K \ra \Gal(L_0/K) \ra \Stab \subset G$. 

Given a continuous
homomorphism $\chi : G_K \ra G$, we have an $\im(\chi)$-extension $L_0$ of $K$
corresponding to the kernel of $\chi$ via Galois theory.  We let
$L=\Ind_{\im(\chi)}^G L_0$. 
It is straightforward to check that these two constructions are inverse to each other.
\end{proof}

\begin{proposition}\label{P:autG}
 A $G$-structured algebra  has exactly $|G|$-automorphisms.
\end{proposition}
\begin{proof}
 Consider a $G$-structured $F$-algebra $L$.
Let $\bar{F}$ be the separable closure of $F$.
There are $|G|$ non-zero morphisms $\phi_i : L \ra \bar{F}$.
Let $S_{|G|}$ be the permutations of these $\phi_i$.
We have $G\subset \Aut_F(L)\sub S_{|G|}$.
An automorphism of $L$ as a $G$-structured $F$-algebra
is an element $\sigma\in \Aut_F(L)$ such that $\sigma$ centralizes $G$.
Clearly all elements of $G$ will satisfy this condition since $G$ is abelian.
Since $G$ acts transitively on the idempotents of $L$, $G$ acts transitively in $S_{|G|}$.  
Thus we can relabel the $\phi_i$ by elements of $G$, and $G$ will act by multiplication on the labels.
So if $\sigma \in S_{|G|}$ centralizes $G$, then
$\sigma$ is translation by an element of $G$, and these are just the automorphisms that come from $G$.
\end{proof}

By class field theory, maps $\chi : G_K \ra G$ are in one-to-one correspondence with maps $\chi : J/K^\times \ra G$. 
Given the correspondence of Lemma~\ref{L:corr}, we can also apply
$C$ to characters $\chi : J/K^\times \ra G$.  
We now view a generator
of tame inertia $y_v$ as an element
of $K_v^\times$, and define $c(\chi_v)$ to be $c_G(\chi_v(y_v))$ for $v$ finite and not dividing $n$. For $v$ infinite or dividing $n$,
let $L_v$ be the $G$-structured  \'{e}tale  $K_v$-algebra corresponding
 to the character $\chi_v$, and define
$c(\chi_v)$ to be $c_v(L_v)$.
We say that 
$$
C(\chi)=\prod_{v \textrm{ place of } K} Nv^{c(\chi_v)}
$$
Just as $\Sigma$ denotes local specifications of $G$-structured  $K_v$-algebras
at the places $v\in S$, we let $\phi$ denote a collection of choices $\phi_v : K_v^\times \ra G$ for
all $v\in S$.  We say that $\phi$ \emph{corresponds to $\Sigma$} if
each $\phi_v$ corresponds to $\Sigma_v$ via Lemma~\ref{L:corr}.

\subsection{Generating functions and Euler products}\label{SS:Euler}
For now, we will assume $C$ is an arbitrary counting function, and in Lemma~\ref{L:fair}, we will first see
how fairness plays a role in our analysis.
Also, for now we will consider one local specification $\Sigma$ (not necessarily viable) on a finite set
$S$ of places of $K$ such that $S$ contains all infinite places, places dividing $n$, and
so that the finite places of $S$ generate the class group of $K$.  
In particular, if $o_S$
is the ring of $S$-integers of $K$ (elements of $K$ with non-negative valuation at all places not in $S$),
then $o_S$ has class number 1.

We define the generating functions
$$
 N_{C,G}(s, \Sigma)= \sum_{\substack{\textrm{$G$-extensions } L/K \\ \forall v\in S \ L_v\isom\Sigma_v }}
\frac{1}{C(L)^{s}} \quad \textrm{ and }\quad N_{C,G}(s, \phi)= \sum_{\substack{\chi : J/K^\times \ra G\\ \forall v\in S \ \chi_v=\phi_v \\ \chi \textrm{ surjective}}}
\frac{1}{C(\chi)^{s}}.
$$
By Lemma~\ref{L:corr},
for $\phi$ corresponding to $\Sigma$ we have $ N_{C,G}(s, \Sigma)= N_{C,G}(s, \phi)$. 
It will be easier to work without the restriction that our characters are surjective, so we define the following generating function:
$$
  F_{C,G}(s, \phi)= \sum_{\substack{\chi : J/K^\times \ra G\\ \forall v\in S \ \chi_v=\phi_v }}
\frac{1}{C(\chi)^{s}}.
$$

For a subgroup $H$ of $G$, we define $C|_H$, a counting function for $H$.  
For $\chi_v:J/K^\times {\ra} H$, we let $c|_H(\chi_v)=c(J/K^\times \stackrel{\chi_v}{\ra} H \subset G)$.
We have that
$$
 F_{C,G}(s, \phi)=\sum_{ H \textrm{ subgroup of } G} N_{C|_H,H}(s, \phi).
$$
We can use M\"{o}bius inversion (as in Wright's work \cite[Section 2]{Wright}) to write
$$
N_{C,G}(s, \phi)= \sum_{ H \textrm{ subgroup of } G} \mu(H,G)  F_{C|_H,H}(s, \phi),
$$
where $\mu(H,G)$ is a constant and $\mu(G,G)=1$.
(This is just solving an upper triangular system of linear equations.)
Thus, by studying the $F_{C,G}$ we can recover information about the $N_{C,G}$.

A character $\chi : J \ra G$ is determined by a collection of
$\chi_v : K_v^\times \ra G$ for all places $v$ of $K$, but not all $\chi$ factor through $J/K^\times$.
However, we can use the following.
\begin{lemma}\label{L:ideleclass}
If $o_S$ has class number $1$, then the natural map
$J_S / o_S^\times \ra J/K^\times$ is an isomorphism.
\end{lemma}
\begin{proof}
 Since $J_S\cap K^\times=o_S^\times$, the map is injective.
Let $x\in J$.  Then, since $o_S$ has class number 1, 
we can find an element of $K$ with specified valuation at all places
outside $S$.  In particular, we can find a $y\in K^\times$ such that
$yx \in J_S$.
\end{proof}

We can then rewrite
$$
  F_{C,G}(s, \phi)= \sum_{\substack{\chi : J_S / o_S^\times \ra G\\ \forall v\in S \ \chi_v=\phi_v }}
\frac{1}{C(\chi)^{s}}.
$$
We shall study characters on $J_S$, and then check their behavior on the finitely generated group $o_S^\times$ to see if they factor
through $J_S / o_S^\times$.  
Let $\mathcal{A}=\prod_{i=1}^k o^\times_S/o^{n_i}_S$.  Given a $\chi : J_S \ra G$, with
projection $\chi_i : J_S \ra \Z/n_i$ (or the same from $K_v^\times$ or $o_v^\times$), and an $\epsilon=(\epsilon_1,\dots,\epsilon_k)\in\mathcal{A}$,
we define $\dot\chi(\epsilon)=\prod_{i=1}^k \zeta_{n_i}^{\chi_i(\epsilon_i)}$,
where we evaluate $\chi_i(\epsilon_i)$ using the natural map $o^\times_S \ra J_S$ (or to $K_v^\times$ or $o_v^\times$).
Note that the map $\chi$ has its image in $G$, the map $\chi_i$ has its
image in $\Z/n_i$, and the map $\dot\chi$ has its image in the complex roots of unity.
We define the twists
$$
  F_{C,G}(s, \epsilon, \phi)= \sum_{\substack{\chi : J_S \ra G\\ \forall v\in S \ \chi_v=\phi_v }}
\frac{\dot{\chi}(\epsilon)}{C(\chi)^{s}},
$$
which we use with the following corollary of Lemma~\ref{L:ideleclass} (motivated by \cite[Equation (3.2)]{Wright}).
\begin{corollary}\label{C:eps}
 We have  $F_{C,G}(s, \phi)=\frac{1}{|\cA|}\sum_{\epsilon \in \cA }  F_{C,G}(s,\epsilon, \phi)$.
\end{corollary}
\begin{proof}
 We rearrange the sum to obtain
$$
\sum_{\epsilon \in \cA }  F_{C,G}(s,\epsilon, \phi)= \sum_{\substack{\chi : J_S \ra G\\ \forall v\in S \ \chi_v=\phi_v }}
\frac{1}{C(\chi)^{s}}
\sum_{\epsilon_1\in o^\times_S/o^{n_1}_S} \zeta_{n_1}^{\chi_1(\epsilon_1)} \cdots \sum_{\epsilon_k\in o^\times_S/o^{n_k}_S}
\zeta_{n_k}^{\chi_k(\epsilon_k)} .
$$
We note that $\zeta_{n_i}^{\chi_i(\epsilon_i)}$ is a complex valued character on the finite group
$ o^\times_S/o^{n_i}_S$ and thus 
$\sum_{\epsilon_i\in o^\times_S/o^{n_i}_S} \zeta_{n_i}^{\chi_i(\epsilon_i)} $
is $|o^\times_S/o^{n_i}_S|$ if $\chi_i$ is trivial on $o_S^\times$ and 0 otherwise.
\end{proof}

The $F_{C,G}(s, \epsilon, \phi)$ are convenient to work with because they have Euler products (as in \cite[Equation (3.4)]{Wright})
$$
F_{C,G}(s, \epsilon, \phi)=\prod_{v\not\in S}\left(\sum_{\chi_v : o_v^\times \ra G} \frac{\dot{\chi_v}(\epsilon)}{Nv^{c(\chi_v)s}} \right) \prod_{v\in S}\frac{\dot{\phi_v}(\epsilon)}{Nv^{c(\phi_v)s}}.
$$
In this paper, all products over $v\not\in S$ are products over the places of $K$ not in $S$.

\subsection{Proof of Main Theorem~\ref{T:prob}}\label{SS:poles}

We will now see how Theorem~\ref{T:prob} will follow from three lemmas, all of which will be proven in 
Section~\ref{S:cont}.
Recall that $m=\min_{g\in G\setminus\{1\}} c_G(g)$ and $\mathfrak{M}=c_G^{-1}(m)$.
 We will prove the following lemma by
relating $F_{C,G}(s, \epsilon, \phi)$ to $L$-functions whose analytic behavior we already know.
\begin{lemma}\label{L:cont}
For any counting function $C$, the product
$F_{C,G}(s, \epsilon, \phi)$
absolutely converges in $\Re(s)> \frac{1}{m}$ and has a meromorphic continuation
to $\Re(s)\geq \frac{1}{m}$, analytic away from $s=\frac{1}{m}$.
The pole of $F_{C,G}(s, 1, \phi)$ at $s=\frac{1}{m}$ is of order
$$
\sum_{g\in \mathfrak{M}} \frac{1}{[K(\zeta_{r_g}):K]},
$$
where $r_g$ is the order of $g$ in $G$.
\end{lemma}
Thus we also obtain a meromorphic continuation to $\Re(s)\geq \frac{1}{m}$ for
$F_{C,G}(s,  \phi)$ and $N_{C,G}(s,  \phi)$.
Lemma~\ref{L:cont} will allow us to use a Tauberian theorem (see \cite[Corollary, p. 121]{N})
to find the probabilities $\Pr_C$.  In the application of the Tauberian theorem,
we will need to know which terms of $F_{C,G}(s,\phi)$ contribute to the main pole, and the following lemma will tell us just that.
\begin{lemma}\label{L:whopoles}
 For a counting function $C$, there is a subgroup $\E$ of $\cA$ such that
if $\epsilon \in \E$ then $F_{C,G}(s, \epsilon, \phi)$ has a pole of the same order at $s=\frac{1}{m}$ as
$F_{C,G}(s, 1, \phi)$, and
if $\epsilon \not\in \E$ then
$F_{C,G}(s, \epsilon, \phi)$ has a pole of smaller order (possibly equal to zero) than that
of
$F_{C,G}(s, 1, \phi)$.
\end{lemma}
The following lemma will allow us to simplify the probabilities we obtain into a reasonable form
for fair counting functions.
\begin{lemma}\label{L:faireps}
 If $C$ is fair, $v\not\in S$, and $\chi_v : o_v^\times \ra G$, then for all $e\in \E$, we have
$\chi_v(e)=0$.
\end{lemma}

Using these lemmas, we can prove Theorem~\ref{T:prob}.

\begin{proof2}{\it{}of Theorem~\ref{T:prob}}
Assume $C$ is fair.
 We have the Euler product
\begin{align*}
 F_{C,G}(s,\epsilon, \phi)=\prod_{v\not\in S}\left(\sum_{\chi_v : o_v^\times \ra G} \frac{\dot{\chi_v}(\epsilon)}{Nv^{c(\chi_v)s}}
 \right) \prod_{v\in S}  \frac{\dot{\phi_v}(\epsilon)}{Nv^{c(\chi_v)s}}. 
\end{align*}
If $e\in\E$, and $\epsilon\in\cA$, Lemma~\ref{L:faireps} implies
\begin{equation}\label{E:firstfair}
 F_{C,G}(s,e\epsilon, \phi)=F_{C,G}(s,\epsilon, \phi)\prod_{v\in S}  \dot{\phi_v}(e).
\end{equation}
Thus,
\begin{align*}
 F_{C,G}(s, \phi)&=\frac{1}{|\cA|}\sum_{\epsilon \in \cA }  F_{C,G}(s,\epsilon, \phi)\\
&=\frac{1}{|\cA|}\sum_{\epsilon \in \cA/\E } \sum_{e\in\E}  F_{C,G}(s,\epsilon , \phi)\prod_{v\in S} 
 \dot{\phi_v}(e)\\
&=\frac{1}{|\cA|}\sum_{\epsilon \in \cA/\E }  F_{C,G}(s,\epsilon , \phi) \sum_{e\in\E}  \prod_{v\in S} 
 \dot{\phi_v}(e),
\end{align*}
where $\cA/\E$ denotes a set of coset representatives for the quotient of $\cA$ by $\E$.

For $g$ and $h$ meromorphic functions on $\Re(s)\geq \frac{1}{m}$, analytic away from $s=\frac{1}{m}$,
 we use $g\sim_m h$ to denote that $g-h$ has a pole at $\frac{1}{m}$
of lesser order then the pole of $g$ (or that at $\frac{1}{m}$ $g-h$ has no pole and $g$ has a pole).

{\bf Case I:}
If $\prod_{v\in S} \dot{\phi_v}$ is not the trivial character on $\E$, we have
$$\sum_{e\in\E}  \prod_{v\in S} 
 \dot{\phi_v}(e)=0$$ and thus
$F_{C,G}(s, \phi)=0$.  This means that there are no $\chi : J/K^\times \ra G$ that for all $v\in S$
have $\chi_v=\phi_v$, and thus $\phi$ is associated to an inviable $\Sigma$.

{\bf Case II:} 
If $\prod_{v\in S} \dot{\phi_v}$ is the trivial character on $\E$.  Then, 
\begin{align*}
 F_{C,G}(s, \phi)
&=\frac{|\E|}{|\cA|}\sum_{\epsilon \in \cA/\E }  F_{C,G}(s,\epsilon , \phi)\\
&\sim_m \frac{|\E|}{|\cA|} F_{C,G}(s,1 , \phi) \quad \textrm{by Lemma~\ref{L:whopoles}.}
\end{align*}
In particular, $F_{C,G}(s, \phi)$ has a pole of order $\sum_{g\in \mathfrak{M}} \frac{1}{[K(\zeta_{r_g}):K]}$ (from Lemma~\ref{L:cont}) at $s=\frac{1}{m}$.

Now we can analyze the pole at $\frac{1}{m}$ of $N_{C,G}(s, \phi)$.  Recall, we can write
$$
N_{C,G}(s, \phi)= \sum_{ H \textrm{ subgroup of } G} \mu(H,G)  F_{C|_H,H}(s, \phi),
$$
By Lemma~\ref{L:cont}, we know that for $H$ a proper subgroup of $G$, the maximum order of
a pole of any $F_{C|_H,H}(s, \epsilon, \phi)$ and thus of any $F_{C|_H,H}(s, \phi)$ 
is $\sum_{g\in \mathfrak{M}\cap H} \frac{1}{[K(\zeta_{r_g}):K]}$.   For fair $C$, this is smaller than the order of the pole of $F_{C,G}(s, \phi)$, and thus
$$
N_{C,G}(s, \phi)\sim_m  F_{C,G}(s, \phi) \sim_m \frac{|\E|}{|\cA|} F_{C,G}(s,1 , \phi).
$$
In particular, $N_{C,G}(s, \phi)$ has a pole at $s=\frac{1}{m}$ and thus
 is not identically zero. 
So there are surjective $\chi : J/K^\times \ra G$ that for all $v\in S$
have $\chi_v=\phi_v$.  So, $\phi$ is associated to a viable $\Sigma$.

If we write $T_{C,G}(s)=\prod_{v\not\in S}\left(\sum_{\chi_v : o_v^\times \ra G} \frac{1}{Nv^{c(\chi_v)s}} \right)$, then
$$ N_{C,G}(s, \phi)\sim_m  \frac{|\E|}{|\cA|} T_{C,G}(s)
 \prod_{v\in S}\frac{1}{Nv^{c(\phi_v)s}} .
$$
Note that $T_{C,G}(s)$ does not depend on $\phi$ and has a pole at $\frac{1}{m}$.  
Let $w$ be the order of the pole of $N_{C,G}(s,\phi)$ (or $T_{C,G}(s)$) at $\frac{1}{m}$.  Let $\Sigma$ be associated to $\phi$.
Let
$$N_{C,G}(\Sigma,X)=\#\{\textrm{$L\in E_G(K)$ $\mid$ $L_v\isom \Sigma_v$ for all $v\in S$ and $C(L)<X$}\}.$$
Then, for viable $\Sigma$, using a Tauberian theorem (as in \cite[Corollary, p. 121]{N}), we obtain
a positive finite limit
$$
\lim_{X \ra \infty} \frac{N_G(\Sigma,X)}{X^{\frac{1}{m}}(\log X)^{w-1}}=
\lim_{s\ra \frac{1}{m}}\left[N_{C,G}(s, \Sigma)(s-\frac{1}{m})^w\right] \frac{m}{\Gamma(w)},
$$
where $\Gamma$ is the Gamma function.  Summing over the finitely many $\Sigma$ on $S$, we have
$$
\lim_{X \ra \infty} \frac{\#\{L\in E_G(K)|C(L)<X\}}{X^{\frac{1}{m}}(\log X)^{p-1}}
$$
is a positive finite  constant.  Thus for viable $\Sigma$ on $S$, we have $\Pr_C(\Sigma)>0$.

It follows that  for a fair counting function $C$ and $\Sigma$ and $\Sigma'$ viable local specifications on $S$, we have
$$
\frac{{\Pr}_C(\Sigma)}{{\Pr}_C(\Sigma')}=
\lim_{X\ra \infty} \frac{N_G(\Sigma_1,X)}{N_G(\Sigma_2,X)}=\frac{\prod_{v\in S} \frac{1 }{Nv^{c(\Sigma_v)/m}}
  }{\prod_{v\in S} \frac{1}{Nv^{c(\Sigma'_v)/m}}}.
$$
We have required that $S$ is sufficiently large
to contain certain places depending on $G$ and $K$ and from our requirements it follows that $S_0\subset S$.
Thus, since $\Sigma$ and $\Sigma'$ are viable, we have $\tilde{\Sigma},\tilde{\Sigma}'\subset A$
and 
${P}_C(\tilde{\Sigma}|A)=\frac{{P}_C(\tilde{\Sigma})}{P_C(A)}$.
We then conclude that $\frac{{\Pr}_C(\Sigma)}{{\Pr}_C(\Sigma')}=\frac{{P}_C(\tilde{\Sigma})}{{P}_C(\tilde{\Sigma'})}=
\frac{{P}_C(\tilde{\Sigma}|A)}{{P}_C(\tilde{\Sigma'}|A)}$.
This proves Theorem~\ref{T:prob} in the case $S$ that is sufficiently large.

Consider a local specification $\Sigma'$ on $S'\subset S$.
Then, we see
$$
{\Pr}_C(\Sigma')=\sum_{\substack{\textrm{viable $\Sigma$ on $S$},\\ \textrm{restricting to $\Sigma'$ on $S'$}}} {\Pr}_C(\Sigma)=
\sum_{\substack{\textrm{viable $\Sigma$ on $S$},\\ \textrm{restricting to $\Sigma'$ on $S'$}}} \frac{{P}_C(\tilde{\Sigma})}{P_C(A)}=
\frac{{P}_C(\tilde{\Sigma'}\cap A)}{P_C(A)},
$$
which proves Theorem~\ref{T:prob}.

\end{proof2}

\subsection{Analytic continuation of $F_{C,G}(s, \epsilon, \phi)$}\label{S:cont}

In this section we prove Lemmas~\ref{L:cont}, \ref{L:whopoles}, and \ref{L:faireps}, the content of which we now remind the reader.

{\it For any counting function $C$, the product
$F_{C,G}(s, \epsilon, \phi)$
absolutely converges in $\Re(s)> \frac{1}{m}$ and has a meromorphic continuation
to $\Re(s)\geq \frac{1}{m}$, analytic away from $s=\frac{1}{m}$.
The pole of $F_{C,G}(s, 1, \phi)$ at $s=\frac{1}{m}$ is of order
$$
\sum_{g\in \mathfrak{M}} \frac{1}{[K(\zeta_{r_g}):K]},
$$
where $r_g$ is the order of $g$ in $G$.

 For a counting function $C$, there is a subgroup $\E$ of $\cA$ such that
if $\epsilon \in \E$ then $F_{C,G}(s, \epsilon, \phi)$ has a pole of the same order at $s=\frac{1}{m}$ as
$F_{C,G}(s, 1, \phi)$, and
if $\epsilon \not\in \E$ then
$F_{C,G}(s, \epsilon, \phi)$ has a pole of lesser (possibly zero) order than
$F_{C,G}(s, 1, \phi)$.

 If $C$ is fair, $v\not\in S$, and $\chi_v : o_v^\times \ra G$, then for all $e\in \E$, we have
$\chi_v(e)=0$.}

We see easily that $F_{C,G}(s, \epsilon, \phi)$ 
(as well as all other products we consider in this subsection)
converges absolutely and uniformly on $\Re(s)>\frac{1}{m}$. 
 So, we will investigate
the behavior at $\frac{1}{m}$ by manipulating the Euler product for $F_{C,G}(s, \epsilon, \phi)$ until it resembles a product of $L$- functions.
This strategy was motivated by the work of Taylor \cite[Section 3]{Taylor}, who related 
$F_{C,G}(s, \epsilon, \phi)$ to $L$-functions for $C$ the conductor and $G$ cyclic, though we face additional
challenges both from general $C$ and $G$ not necessarily cyclic.

  We use the following lemma to interchange sums and 
products, which is possible because we are only looking for  behavior at $\frac{1}{m}$ and so higher order terms will not contribute.
For $g$ and $h$ analytic functions on $\Re(s)> \frac{1}{m}$, we use $g\approx_m h$ to
denote that 
$\frac{g}{h}$ has an analytic continuation to $\Re(s)\geq \frac{1}{m}$.
\begin{lemma}\label{L:Pv}
Let $m$ and $M$ be positive reals.
Let $K$ be a number field, for each place $v$ of $K$, let
$P_v(x)=1+\sum_{i=1}^{M} b_{v,i}x^{\alpha_{v,i}}$, where $m\leq \alpha_{v,i}$ and
$b_{v,i} \in \C$ with $|b_{v,i}|\leq M$.  Then
for some large $Y$ we have
 $$\prod_v P_v(Nv^{-s}) \approx_m \prod_{\substack{v\\ \Nv > Y}}\prod_{i=1}^{M}
\left(1+ b_{v,i}\Nv^{-\alpha_{v,i}s}\right)$$ (where the
products over $v$ are over all finite places $v$ of $K$ satisfying the condition).
\end{lemma}

\begin{proof}
We can bound the absolute value of each factor of $\prod_v P_v(Nv^{-s})$
by $ 1 + M^2 \Nv^{-ms}$ and each factor of $\prod_v \prod_{i=1}^{M}
\left(1+ b_{v,i}\Nv^{-\alpha_{v,i}s}\right)$ by $(1+M\Nv^{-ms})^M$, and thus
both products converge absolutely on $\Re(s)>\frac{1}{m}$.
For sufficiently large $v$, the function $\prod_{i=1}^{M}
\left(1+ b_{v,i}\Nv^{-\alpha_{v,i}s}\right)$ has absolute value at least $\frac{1}{2}$ everywhere on  $\Re(s)\geq \frac{1}{m}$.
For those $v$,
$$
\left|\frac{P_v(Nv^{-s})}{\prod_{i=1}^{M}
\left(1+ b_{v,i}\Nv^{-\alpha_{v,i}s}\right)}\right|\leq 1+ \frac{2^M M^M \Nv^{-2ms}}{|\prod_{i=1}^{M}
\left(1+ b_{v,i}\Nv^{-\alpha_{v,i}s}\right)|}
\leq 1+ 2^{M+1} M^M \Nv^{-2ms}.
$$
Thus we conclude the lemma.
\end{proof}

Now, we set our notation for the rest of the proof of Lemmas~\ref{L:cont}, \ref{L:whopoles}, and \ref{L:faireps}.

\begin{notation}
A \emph{division}
of $G$ is a set of all the invertible multiples of some element $x\in G$, in other words 
 $\{y \mid y=ex \textrm{ and } x=fy \textrm{ for some } e,f\in\Z \}$.
Let $\Div(G)$ be the set of non-identity divisions of $G$.    For an element $g\in G$, let $r_g$ be its order
and for $d\in\Div(G)$, let $r_d$ be the order of any element of $d$.
Recall that any map from $o_v^\times$ to a finite group of order relatively prime to $v$ factors through
$(o_v/v)^\times$.

We now make a specific choice, for all places $v\nmid |G|$, of a 
generator $y_v$ of the tame inertia group of $K_v$  (which is isomorphic to $(o_v/v)^\times$).
Our choice is that $y_v \equiv \zeta_{\Nv-1} \pmod{v}$, where $\zeta_{\Nv-1}$
is the in the primitive $({\Nv-1})$th root of unity we fixed just before Section~\ref{SS:G}.

Since $c(\chi_v)$ only depends on the division of $\chi_v(y_v)$, for a division $d$
we can write $c(d)$ to denote $c(\chi_v)$ for any $\chi_v$ that sends $y_v$ to an element of $d$.

\end{notation}

We now rearrange $F_{C,G}(s,\epsilon, \phi)$ as follows
\begin{align*}
F_G(s,\epsilon, \Sigma)&\approx_m \prod_{v\not\in S} \left( \sum_{\chi_v : o_v^\times \ra G} \frac{\dot{\chi_v}(\epsilon)}{Nv^{c(\chi_v)s}} \right)\\
& = \prod_{v\not\in S} \left( 1+ \sum_{d\in\Div(G)}  \sum_{g\in d} \sum_{\substack{\chi_v : o_v^\times \ra G\\ 
\chi_v(y_v)=g }} 
\frac{\dot{\chi_v}(\epsilon)}{Nv^{c({d})s}} \right).
\end{align*}\
The sum over $\chi_v : o_v^\times \ra G$ such that 
$\chi_v(y_v)=g$ has at most one term, but we keep the summation sign for notational convenience.
So we have
\begin{align*}
F_G(s,\epsilon, \Sigma)
&\approx_m \prod_{\substack{v\not\in S\\ Nv> Y }}  \prod_{d\in\Div(G)} \prod_{g\in d} \left( 1+  \sum_{\substack{\chi_v : o_v^\times \ra G\\ \chi_v(y_v)=g}} 
\frac{\dot{\chi_v}(\epsilon)}{Nv^{c({d})s}} \right) \quad \mbox{by Lemma~\ref{L:Pv}}\\
&=   \prod_{d\in\Div(G)} \prod_{\substack{v\not\in S\\ \Nv \equiv 1 \pmod{r_d}\\Nv >Y }} \prod_{g\in d}\left( 1+  \frac{1}{Nv^{c({d})s}}\sum_{\substack{\chi_v : o_v^\times \ra G\\ 
\chi_v(y_v)=g}} 
\dot{\chi_v}(\epsilon) \right).
\end{align*}
Only $v$ with
$\Nv \equiv 1 \pmod{r_d}$ have $\chi : o_v^\times \ra G$ such that
$\chi_v(y_v)\in d$.

Now we prove the following lemmas in order to evaluate the term
$
\dot{\chi_v}(\epsilon)$ in the above.  Our strategy to evaluate 
$\dot{\chi_v}(\epsilon)$ is motivated by the work of Taylor \cite{Taylor}, who calculated the order
of $\dot{\chi_v}(\epsilon)$ for $G$ cyclic.  For non-cyclic $G$, we need to take advantage of our choice of $y_v$.
\begin{lemma}\label{L:whyy}
 We have 
$$
\zeta_{\Nv-1}=\frac{\Frob_v(y_v^{1/(\Nv-1)})}{y_v^{1/(\Nv-1)}},
$$
where the Frobenius is in the Galois group of the maximal unramified extension of $K_v$. 
\end{lemma}
\begin{proof}
 Note
that $K_v$ contains the $(\Nv-1)$th roots of unity and so $\frac{\Frob_v(y_v^{1/(\Nv-1)})}{y_v^{1/(\Nv-1)}}$ 
does not depend on the choice of root of $y_v$.
We know that both $\zeta_{\Nv-1}$ and $\frac{\Frob_v(y_v^{1/(\Nv-1)})}{y_v^{1/(\Nv-1)}}$
are $(\Nv-1)$th roots of unity, and that those roots of unity inject into $(o_v/v)^\times$.
Thus we can prove the lemma modulo $v$. There we have
$$
\frac{\Frob_v(y_v^{1/(\Nv-1)})}{y_v^{1/(\Nv-1)}}\equiv y_v^{(\Nv-1)/(\Nv-1)}\equiv y_v \equiv \zeta_{Nv-1},
$$
where the last equality is by choice of $y_v$.
\end{proof}

\begin{lemma}\label{L:findfrob}
Let $v \nmid n\infty$ and $\chi_v(y_v)=g$.
Suppose the projections of $g$ 
to the $\Z/n_i\Z$ are $\frac{n_i k_i}{\ell_i} \in \Z/n_i\Z$, where $\ell_i \mid n_i$ and $(k_i,\ell_i)=1$.
Let $\epsilon^g$ be notation for $\prod_{i=1}^k \epsilon_i^{k_i/\ell_i}$,
and let $w_v$ be a prime of $K(\zeta_{r_g})$ over $v$.
 Then
$$\dot{\chi_v}(\epsilon)=\prod_{i=1}^k \frac{\Frob_{w_v}(\epsilon_i^{k_i/\ell_i})}{\epsilon_i^{k_i/\ell_i}}
= \frac{\Frob_{w_v}(\epsilon^g)}{ \epsilon^g},$$
where the Frobenius is in the Galois group of the maximal extension of $K(\zeta_{r_g})$ unramified outside $S$.
\end{lemma}
\begin{proof}
Note that $\ell_i\mid r_g$.  
Since $\chi_v$ factors through $(o_v/v)^\times$ and $y_v$ has order $Nv-1$ in $(o_v/v)^\times$, we also
have that $r_g \mid Nv-1$.
In $(o_v/v)^\times$, write $\epsilon_i=y_v^{b_i}$ so that in $K_v$ we have $\epsilon_i u_i=y_v^{b_i}$,
where $u_i$ is a unit congruent to 1 modulo $v$.
We have that 
\begin{align*}
 \dot{{\chi_v}_i}(\epsilon_i)=\zeta_{n_i}^{{\chi_v}_i(\epsilon_i)} 
=\zeta_{n_i}^{b_i {\chi_v}_i(y_v)} =\zeta_{n_i}^{ \frac{n_i  k_i b_i}{\ell_i}}
=\zeta_{\ell_i}^{   k_i b_i} 
= \zeta_{\Nv-1}^{ \frac{ (\Nv-1)k_i b_i }{\ell_i}}.
\end{align*}
From Lemma~\ref{L:whyy}, we have $\zeta_{\Nv-1}=\frac{\Frob_v(y_v^{1/(\Nv-1)})}{y_v^{1/(\Nv-1)}}$,
where the Frobenius is in the Galois group of the maximal unramified extension of $K_v$.  
Thus
$$
 \dot{{\chi_v}_i}(\epsilon_i)=\left(\frac{\Frob_v(y_v^{1/(\Nv-1)})}{y_v^{1/(\Nv-1)}}\right)^{\frac{ (\Nv-1)k_i b_i }{\ell_i}}=
\frac{\Frob_v(y_v^{k_i b_i/\ell_i})}{y_v^{k_i b_i/\ell_i}}=
\frac{\Frob_v(\epsilon_i^{k_i/\ell_i})\Frob_v(u_i^{k_i/\ell_i})}{\epsilon_i^{k_i/\ell_i}u_i^{k_i/\ell_i}},
$$
where the Frobenius is still in the Galois group of the maximal unramified extension of $K_v$.
Since $u_i$ is a unit congruent to 1 modulo $v$ and $\ell_i \mid Nv-1$, we have that all the $\ell_i$th roots of $u_i$ are
in $K_v=K(\zeta_r)_{w}$ and that $\Frob_v(u_i^{k_i/\ell_i})=u_i^{k_i/\ell_i}$.
Note that $K_v=K(\zeta_{r_g})_{w_v}$ since $r_g | Nv-1$, and thus we 
can replace $\Frob_v$ with the Frobenius of $w_v$ in $K(\zeta_{r_g})_{w_v}$.
We thus have
$$
 \dot{{\chi_v}_i}(\epsilon_i)=
\frac{\Frob_{w_v}(\epsilon_i^{k_i/\ell_i})}{\epsilon_i^{k_i/\ell_i}}.
$$
Since the $\ell_i$th roots of $\epsilon_i$ are in the
the maximal extension of $K(\zeta_{r_g})$ unramified outside $S$, we can interpret the Frobenius as 
the Frobenius of $w_v$ in the Galois group of the maximal extension of $K(\zeta_{r_g})$ unramified outside $S$
in the 
statement of the Lemma. Note
that $K(\zeta_{r_g})$ contains the $\ell_i$th roots of unity and so $\frac{\Frob_v(\epsilon_i^{k_i/\ell_i})}{\epsilon_i^{k_i/\ell_i}}$ 
does not depend on the choice of root of $\epsilon_i$.
\end{proof}

Using Lemma~\ref{L:findfrob} and its definitions of $\epsilon^g$, $w_v$, and $\Frob$, we have 
\begin{align*}
F_{C,G}(s,\epsilon, \phi)&\approx_m \prod_{d\in\Div(G)} \prod_{\substack{v\not\in S\\ \Nv \equiv 1 \pmod{r_d}\\Nv >Y }} \prod_{g\in d}\left( 1+  \frac{1}{Nv^{c({d})s}}\sum_{\substack{\chi_v : o_v^\times \ra G\\ 
\chi_v(y_v)=g}} 
\dot{\chi_v}(\epsilon) \right)\\
&= \prod_{d\in\Div(G)} \prod_{\substack{v\not\in S\\ \Nv \equiv 1 \pmod{r_d}\\Nv >Y }} \prod_{g\in d} \left( 1+  \frac{1}{Nv^{c({d})s}}
\frac{\Frob_{w_v}(\epsilon^g)}{\epsilon^g} \right).
\end{align*}

We now partition $\Div(G)$ into $\Div^0(\epsilon,G)$,
the divisions whose elements $g$ have $\epsilon^g \in K(\zeta_{r_g})$, and $\Div^+(\epsilon, G)$,
the divisions whose elements $g$ have $\epsilon^g \not\in K(\zeta_{r_g})$.
Let $t(r):=[K(\zeta_r):K]$.
We factor
the last product above into two factors $A(s)$ and $B(s)$, defined below.

We have 
\begin{align*}
A(s)&:=\prod_{d\in\Div^0(\epsilon, G)} \prod_{\substack{v\not\in S\\ \Nv \equiv 1 \pmod{r_d}\\Nv >Y }} \prod_{g\in d} \left( 1+  \frac{1}{Nv^{c({d})s}}
\frac{\Frob_{w_v}(\epsilon^g)}{\epsilon^g} \right)\\ 
&=\prod_{d\in\Div^0(\epsilon, G)} \prod_{\substack{v\not\in S\\ \Nv \equiv 1 \pmod{r_d}\\Nv >Y }} \prod_{g\in d} \left( 1+  \frac{1}{Nv^{c({d})s}}
\right)\\ 
&=\prod_{d\in\Div^0(\epsilon, G)} \prod_{\substack{v\not\in S\\ \Nv \equiv 1 \pmod{r_d}\\Nv >Y }} \prod_{w | v } \left( 1+  \frac{1}{Nv^{c({d})s}}
\right)^{\phi(r_d)/t(r_d)},\\ 
\end{align*}
where the last product is over the primes $w$ of $K(\zeta_{r_d})$ over $v$.
Note that $\frac{\phi(r_d)}{t(r_d)}$ is an integer.
By the standard argument about only degree one primes contributing to the pole, we have
$$
\prod_{\substack{v\not\in S\\ \Nv \equiv 1 \pmod{r_d}\\Nv >Y }} \prod_{w | v } \left( 1+  \frac{1}{Nv^{c({d})s}}
\right) \approx_m \zeta_{K(\zeta_{r_d})}(c(d)s).
$$
Thus
$$
A(s)\approx_m \prod_{d\in\Div^0(\epsilon, G)} \zeta_{K(\zeta_{r_d})}(c(d)s)^{\phi(r_d)/t(r_d)}.
$$

We define 
$$
B(s):=\prod_{d\in\Div^+(\epsilon, G)} \prod_{\substack{v\not\in S\\ \Nv \equiv 1 \pmod{r_d}\\Nv >Y }} \prod_{g\in d} \left( 1+  \frac{1}{Nv^{c({d})s}} 
\frac{\Frob_{w_v}(\epsilon^g)}{\epsilon^g} \right).
$$
Let $N$ be the least common multiple of the $n_i$, and note that since $r_d\mid N$, we have that
$t(r_d)\mid t(N)$.  
We now have
\begin{align*}
B(s)^{t(N)} =\prod_{d\in\Div^+(\epsilon, G)} \prod_{\substack{v\not\in S\\ \Nv \equiv 1 \pmod{r_d}\\Nv >Y }} \prod_{g\in d}
\prod_{w | v }
 \left( 1+  \frac{1}{Nv^{c({d})s}}
\frac{\Frob_{w}(\epsilon^g)}{\epsilon^g} \right)^{t(N)/t(r_d)},
\end{align*}
where the last product is over the primes $w$ of $K(\zeta_{r_d})$ over $v$.
For $d\in \Div^+(\epsilon, G)$ we have that  $K(\zeta_{r_d},\epsilon)/K(\zeta_{r_d})$ is abelian and non-trivial.
Thus 
there is a non-trivial Hecke character $\theta_{\epsilon^g}$ for $K(\zeta_{r_d})$ 
such that $\frac{\Frob_{w}(\epsilon^g)}{\epsilon^g}$ is 
 $\theta_{\epsilon^g}(w)$.
Again by standard arguments we have
$$
\prod_{\substack{v\not\in S\\ \Nv \equiv 1 \pmod{r_d}\\Nv >Y }} \prod_{w | v } \left( 1+  \frac{\Frob_{w}(\epsilon^g)}{\epsilon^g} \frac{1}{Nv^{c({d})s}}
\right)\approx_m L(c(d)s, \theta_{\epsilon^g})
$$
and thus we can write
$$
B(s)^{t(N)} = g(s) \prod_{d\in\Div^+(\epsilon, G)} L(c(d)s, \theta_{\epsilon^g})^{t(N)/t(r_d)},
$$
where $g(s)$ is analytic in $\Re(s)\geq \frac{1}{m}$.
We know that $ L(c(d)s, \theta_{\epsilon^g})$ not
only has an analytic continuation to $\Re(s) \geq \frac{1}{m}$ but is also non-zero in that region.
We can check that $g(s)$ is also non-zero in $\Re(s)\geq \frac{1}{m}$.
Thus $B(s)$ has an analytic continuation to $\Re(s) \geq \frac{1}{m}$.

Thus, we conclude that
\begin{align*}
F_{C,G}(s,\epsilon, \phi)&\approx_m  \prod_{d\in\Div^0(\epsilon, G)} \zeta_{K(\zeta_{r_d})}(c(d)s)^{\phi(r)/t(r_d)}.
\end{align*}

So, $F_{C,G}(s,\epsilon, \phi)$ has a meromorphic continuation to $\Re(s)\geq\frac{1}{m}$ analytic away from
$s=\frac{1}{m}$.  The pole of $F_{C,G}(s,\epsilon, \phi)$ at $\frac{1}{m}$ is of order
$$
\sum_{\substack{d\in\Div^0(\epsilon, G)\\ c(d)=m}} \frac{\phi(r_d)}{t(r_d)}=\sum_{\substack{d\in\Div^0(\epsilon, G)\cap \mathfrak{M}\\}} \frac{\phi(r_d)}{t(r_d)}
=\sum_{\substack{d\in\Div^0(\epsilon, G)\cap \mathfrak{M}\\}} \sum_{g \in d}\frac{1}{t(r_g)}
=\sum_{g \in G(\epsilon)\cap \mathfrak{M} }\frac{1}{t(r_g)},
$$
 where $G(\epsilon)$ is the set of $g\in G$ such that $\epsilon^g \in K(\zeta_{r_g})$.  Note that 
$G(1)=G$.  
This proves Lemma~\ref{L:cont}.
Clearly the maximal order pole among terms
 $F_{C,G}(s,\epsilon, \phi)$ is in $F_{C,G}(s,1, \phi)$, and any other
$F_{C,G}(s,\epsilon, \phi)$ has that same order pole if and only if $\mathfrak{M}\subset G(\epsilon)$.  
Let $\E$ be the elements $\epsilon \in \cA$ such that $\mathfrak{M}\subset G(\epsilon)$.
It is easy to see $\E$ is a subgroup, and this proves Lemma~\ref{L:whopoles}.
Lemma~\ref{L:faireps} will follow from the next result.

\begin{lemma}\label{L:fair}
 For a fair counting function $C$, and $\epsilon\in \E$, we have 
${\epsilon_j}^{1/r}\in K(\zeta_r)$ for all $r\mid n_j$.  
\end{lemma}

\begin{proof}
 Fix a $j$ with $1\leq j \leq k$ and an $r$ dividing $n_j$.  
Let $g$ be the element of $G$ with $j$th projection $\frac{n_j}{r}$ and all other projections 0.
Since $g$ is of order $r$ and $C$ is fair, we can write
$\sum_{s=1}^\ell  g_s=g$, where $g_s$ are elements of $\mathfrak{M}$
and all $g_s$ have order dividing $r$.  Write $g_s=(g_{s,1},\dots,g_{s,k})$ according to our chosen factorization of $G$.
We can write $g_{s,i}=\frac{n_i h_{s,i}}{\ell_{s,i}}$ with $(h_{s,i},\ell_{s,i})=1$.  Since 
$g_s$ is of order dividing $r$, we must have $\ell_{s,i}|r$. 
Thus by definition of $\E$ we have $$\epsilon^{g_s}=\prod_{i=1}^k 
{\epsilon_i^{h_{s,i}/\ell_{s,i}}}\in K(\zeta_{r}).$$ 
We then see that
$$\prod_{s=1}^\ell \prod_{i=1}^k 
\epsilon_i^{h_{s,i}/\ell_{s,i}}\in K(\zeta_{r}).$$
By the choice of the $g_s$, we have that
$\sum_{s=1}^\ell \frac{n_i h_{s,i}}{\ell_{s,i}}$ (as a sum in $\Z/n_i$) is
$\frac{n_j}{r}$ if $i=j$ and 0 otherwise.
Equivalently,
$\sum_{s=1}^\ell  \frac{h_{s,i}}{\ell_{s,i}}$ (as a sum in $\Q/\Z$) is
$\frac{1}{r}$ if $i=j$ and 0 otherwise.
Thus, we conclude that $\prod_{s=1}^\ell \prod_{i=1}^k 
\epsilon_i^{h_{s,i}/d_{s,i}}$ is  ${\epsilon_j}^{1/r}$
times an element of $K^\times$, and thus ${\epsilon_j}^{1/r}\in K(\zeta_{r})$.
\end{proof}

Suppose $C$ is fair, $v\not\in S$, and  we have a $\chi:o_v^\times \ra G$ of order $r$,
with projection to $\Z/n_i\Z$ of order $\ell_i$.  Then
$\Nv \equiv 1 \pmod{r}$, and thus for all $i$ we have  $K_v=K_v(\zeta_{\ell_i})$.
So, for all $\epsilon \in\E$, we have ${\epsilon_j}^{1/\ell_j}\in K(\zeta_{\ell_j})$, which implies that
${\epsilon_j}^{1/\ell_j}\in K_v(\zeta_{\ell_j})=K_v$, and thus $\epsilon_j$ is a $\ell_j$th power
in $o_v^\times$ for all $j$.  We conclude $\dot{\chi_v}(\epsilon)=0$, which proves Lemma~\ref{L:faireps}.

\begin{remark}\label{ERemark}
 By definition, $\E$ depends on our choice of $C$.  However, given that $C$ is fair,
by Lemma~\ref{L:fair}, we see that for $\epsilon\in\E$ we have
$\epsilon^g \in K(\zeta_{r_g})$ for \emph{all} $g\in G$.  If
 $\epsilon\in\cA$ is such that $\epsilon^g 
\in K(\zeta_{r_g})$ for \emph{all} $g\in G$, then $\epsilon\in \E$.
Thus if $C$ is fair, we see that $\E$ is the subgroup of $\epsilon$ such that
$\epsilon^g 
\in K(\zeta_{r_g})$ for all $g\in G$, and thus does not depend on $C$.
\end{remark}

\section{Counting by conductor}\label{S:constant}
In the proof of Theorem~\ref{T:prob} in Section~\ref{SS:poles}, we have implicitly found the asymptotics of 
$$
N_{C,G}(X):=\#\{L\in E_G(K)|C(L)<X\},
$$
for any fair counting function $C$.  We collect that result here.
Recall the definition of $S_0$ from Section~\ref{S:anal} as follows.
Let $\eta_i=\zeta_{2^i}+\zeta_{2^i}^{-1}$, where $\zeta_{2^i}$ is a primitive $2^i$th root of unity.
Let $s_K$ be maximal such that $\eta_{s_K}\in K$.  If $2^{{s_K}+1}$ does not divide the exponent of $G$,
then let $S_0(K)=\emptyset$.  Otherwise, let $S_0(K)$
be the set of primes $\wp$ of $K$ dividing 2 such that none of $-1$, $2+\eta_{s_K}$ and $-2-\eta_{s_K}$
are squares in $K_\wp$.

\begin{theorem}\label{T:constant}
For a fair counting function $C$, we have
\begin{align*}
 &\lim_{X \ra \infty} \frac{N_{C,G}(X)}{X^{m_C}(\log X)^{w_{K,C}-1}}\\
&=\frac{\Sp(K,G)}{m_C^{w_{K,C}-1}(w_{K,C}-1)!|G|^{|S_0(K)|}
\prod_i{|o_K^\times/o_K^{n_i}|}}
\prod_{\substack{v\not\in S_0(K)\\ v \textrm{ finite}}}\left( \left(\sum_{\chi_v : o_v^\times \ra G} \frac{1}{\Nv^{c(\chi_v)/m_C}}\right)\left(1-\frac{1}{\Nv}\right)^{w_{K,C}}\right)\\
&\qquad\qquad\qquad\cdot
\left( \sum_{\substack{\Sigma \text{ viable local}\\\textrm{spec. of $G$-structured} \\ \textrm{algebras on } S_0(K)}} \prod_{v\in S_0(K)} \frac{1}{\Nv^{c(\Sigma_v)/m_C}}\left(1-\frac{1}{\Nv}\right)^{w_{K,C}}
\right)
\prod_{v\mid\infty} (\sum_{G_{K_v} \ra G } 1 ),
\end{align*}
where
\begin{align*}
 G&=\Z/n_1\Z\times\dots\times\Z/n_k\Z,\\
m_C &= \min_{g\in G\setminus\{0\}}{c_G(g)},\\
\mathfrak{M}&=c_G^{-1}(m),\\
r_g &\textrm{ is the order of } g\in G,\\
\zeta_j &\textrm{ are the $j$th roots of unity},\\
w_{K,C}&=\sum_{g\in \mathfrak{M}} \frac{1}{[K(\zeta_{r_g}):K]},\\
h_{G,K} &\textrm{ is the number of $i$ such that $2^{s_K+1}|n_i$}\\
Sp(K,G)&=\begin{cases}
2^{h_{G,K}} \quad &\textrm{if none of $-1$, $2+\eta_s$ and $-2-\eta_s$ are squares in $K$,}
 \\
1 \quad  \quad &\textrm{otherwise},
\end{cases}\\
o_K &\textrm{ is the ring of integers in $K$},\\
G_F &\textrm{ is the absolute Galois group of $F$},\\
o_v &\textrm{ is the ring of integers of $K_v$},
\end{align*}
and all products are over places of $K$.
\end{theorem}

We can also specialize to the case that the counting function is $\f$, the norm of the conductor to $\Q$.
In this case $m_{\f}=1$ and $\mathfrak{M}=G\setminus\{0\}$ and so the expression in Theorem~\ref{T:constant} simplifies slightly.

\begin{proof}
This result follows from the analysis of Secton~\ref{SS:poles}.  We simplify the constant that one
obtains using that analysis by
applying $$\lim_{s\ra \frac{1}{m}}\left[\zeta_K(sm)(s-\frac{1}{m})\right] =\frac{1}{m}
$$ and the following two lemmas.

 \begin{lemma}
For fair $C$, we have $|\E|=Sp(K,G)$.
 \end{lemma}
\begin{proof}
We know from Lemma~\ref{L:fair} that $\epsilon \in \E$ implies that for all $\ell_i \mid n_i$ we have
$\epsilon_i^{\frac{1}{\ell_i}}\in K(\zeta_{\ell_i})$.  If $v\not\in S$,
and $(\Nv-1, n_i)=\ell_i$, then $\epsilon_i$ is an $n_i$th power in $K_v$ if and only if
it is an $\ell_i$th power in $K_v$.  Since  $\ell_i \mid Nv-1$, we have
$K_v(\zeta_{\ell_i})=K_v$, and thus $\epsilon_i^{\frac{1}{n_i}}\in K_v$.  By \cite[Chapter 10, Theorem 1]{AT},
it follows that either 1) $\epsilon_i$ is an $n_i$th power in $K$ or 2) $\epsilon_i\in b_0^{n_i/2} K^{n_i}$, where
$b_0=2+\eta_{s_K}=(1+\zeta_{2^{s_K}})^{2}$.  Also, the second case only occurs when  none of $-1$, $2+\eta_{s_K}$ and $-2-\eta_{s_K}$
are squares in $K$ and when $2^{s_K+1} \mid n_i$.

Next, we will see that any $\epsilon$ such that $\epsilon_i=b_0^{n_i/2}$ and $2^{s_K+1} \mid n_i$ for
some $i\in I$ and $\epsilon_j=1$ for all $j\not\in I$ is in $\E$.
First note that $b_0$ is a unit at all places not dividing 2, and so
it will be in $o_S^\times$ as long as $S$ contains 2 (which we have required when
$|G|$ is even).
We can reduce to the case that $I=\{i\}$.  Then we need to conclude that
$b_0^{n_i/2}\in K(\zeta_{\ell_i})^{\ell_i}$ for all $\ell_i\mid n_i$.
We can easily reduce to the case that $n_i$ is a power of 2 (for example by choosing the
$n_i$ to be prime powers originally).   We know that 
 $b_0^{n_i/2}\in K^{n_i/2}$.
If $\ell_i=n_i$, then $K(\zeta_{\ell_i})$ contains
$\zeta_{2^{s_K}}$ (because $2^{s_K}\mid n_i$)
and thus $b_0^{n_i/2}$ is a $\ell_i$ th power in $K(\zeta_{\ell_i})$.

We see that $\E$ 
is trivial when any of $-1$, $2+\eta_{s_K}$ and $-2-\eta_{s_K}$
are squares in $K$.
Otherwise, $\E$
contains exactly the $\epsilon$ that have $\epsilon_i=1$
where $2^{s_K+1} \nmid n_i$ and that have $\epsilon_i=1$ or $b_0^{n_i/2}$
at all other $i$.  
From \cite[Chapter 10, Theorem 1]{AT} we know that $b_0^{n_i/2}$
is not an $n_i$th power in $K$ when none of $-1$, $2+\eta_{s_K}$ and $-2-\eta_{s_K}$
are squares in $K$ and $2^{s_K+1} \mid n_i$. 
This proves the lemma.
\end{proof}

The next lemma follows from the fact that a local specification of $G$-structured algebras 
on $S$ containing $S_0(K)$
is viable if and only if its restriction to $S_0$ is viable (see \cite[Chapter 10, Theorem 5]{AT}).
Also recall Lemma~\ref{L:corr}, which gives the correspondence between $G$-structured algebras and 
Galois representations.
 \begin{lemma}
For $S$ containing $S_0(K)$,
\begin{align*}
 &\sum_{\substack{\Sigma \text{ viable local}\\\textrm{spec. of $G$-structured} \\ \textrm{algebras on } S}} \prod_{v\in S} \frac{1}{\Nv^{c(\Sigma_v)/m_C}}=\\
&\prod_{v\in S\setminus S_O} \left(|G|\sum_{\chi_v : o_v^\times \ra G} \frac{1}{\Nv^{c(\chi_v)/m_C}} \right)
\sum_{\substack{\Sigma \text{ viable local}\\\textrm{spec. of $G$-structured} \\ \textrm{algebras on } S_0(K)}} \prod_{v\in S_0(K)} \frac{1}{\Nv^{c(\Sigma_v)/m_C}}.
\end{align*}
 \end{lemma}

\end{proof}

\section{Discriminant Probabilities}\label{S:disc}

For this section, we work with base field $K=\Q$.  We show that when one replaces the conductor
by the discriminant when defining probabilities in Equation~\eqref{E:probdef} (to define what we call
 \emph{discriminant probabilities}), we do not in general have analogs of the nice behavior of Corollary~\ref{C:QCheb}
and Theorem~\ref{T:Qind}.  
When $G$ has prime exponent, the discriminant is  a fixed power of the conductor,
and so we do have analogs of Corollary~\ref{C:QCheb}
and Theorem~\ref{T:Qind}.  
However, in the simplest case when $G$ does not have prime exponent, that is $G=\Z/p^2\Z$ for $p$ prime,
 we find examples of dependence of behaviors at different places (Proposition~\ref{P:notind}), and examples where we do not have the
 Chebotarev probabilities
for unramified splitting behavior (Proposition~\ref{P:notCheb}).
As discussed in the introduction, Wright \cite{Wright} observed
that for $G=\Z/4\Z$ the ratios of probabilities of local behaviors are apparently very complicated.
Our propositions give concrete evidence for the suggestion of Wright that the discriminant probabilities are not well-behaved.
We compute the probabilities for the propositions below in a similar fashion to our work in Section~\ref{S:anal}.

\begin{proposition}\label{P:notind}
Let $p$, $q_1,$ and $q_2$ be primes with $q_i\equiv 1\pmod{p^2}$ for $i=1,2$.
 If $G=\Z/p^2\Z$, then $q_1$ ramifying and $q_2$ ramifying in a random $G$-extension
are not discriminant independent. 
\end{proposition}

\begin{proof}
From Lemma~\ref{L:ideleclass} with $S$ empty, we have $J/\Q^\times\isom\prod_p \Z^\times_p \times \R/\{\pm 1\},$
where the product is over finite places of $\Q$.  We also have the following.
\begin{lemma}\label{L:Q}
 The natural map
$$
\Hom(\prod_p \Z^\times_p \times \R/\{\pm 1\},G)\ra\Hom(\prod_p \Z^\times_p,G),
$$
sending $\chi \mapsto \chi(-,1)$, is an isomorphism.
\end{lemma}

We work as in Section~\ref{S:anal}, but now we let $\Phi$ be a
set
of isomorphism classes of $G$-structured algebras at each place of $S$ instead of just considering a single
$G$-structured algebra. 
For the following computations, we let $G$ be either $\Z/p^2\Z$ or $p\Z/p^2\Z$, and
let $D$ on $\Z/p^2\Z$ be given by $D(L)=|\Disc_\Q L|$, and $D$ on $p\Z/p^2\Z$ be given by
 $D|_{p\Z/p^2\Z}(L)=|\Disc_\Q L|^p$.
In both cases, $m=p(p-1)$.
Let $S$ be a finite set of finite places and let $\Phi$ specify that a
character is unramified at all places in $S$. 
We consider 
$$F_{D,G}(s,\Phi)
:=\sum_{\substack{\chi : J/\Q^\times \ra G\\ \forall v\in S \ \chi_v\in\Phi_v }}
\frac{1}{D(\chi)^{s}}
=\prod_{\ell\not\in S}\left(\sum_{\chi_\ell : \Z_\ell^\times \ra G} \frac{1}{D(\chi_\ell)^{s}} \right) ,$$
where the product is over finite rational primes $\ell$.  We can express $F_{D,G}(s,\Phi)$ as this Euler product by Lemma~\ref{L:Q},
which allows us to count characters from $J/\Q^\times$ by counting characters from $\prod_\ell \Z^\times_\ell$.
We know that $D$ only depends on the restriction of local characters to $\Z^\times_\ell$.
We see that $F_{D,G}(s,\Phi)$ only differs by finitely many factors from the $F_{D,G}(s,1,\phi)$ of Section~\ref{S:anal}
(for any choice of $\phi$), and that
 $F_{D,G}(s,\Phi)/F_{D,G}(s,1,\phi)$ is entire. 
We conclude from Lemma~\ref{L:cont} that
$F_{D,G}(s,\Phi)$ has a pole at $\frac{1}{p(p-1)}$ of order 1, but otherwise can be analytically continued
to $\Re(s)\geq \frac{1}{p(p-1)}$.
As in the end of Section~\ref{S:anal}, we can use a Tauberian theorem to calculate the coefficient sums
$$
F_{D,G}(\Phi,X)=\#\{\textrm{$\chi : J/\Q^\times \ra G$ $\mid$ $\chi_v\in\Phi_v$ for all $v\in S$ and $D(\chi)<X$}\}.
$$

If we let $\Phi^{(q_i)}$ specify that a character is unramified at $q_i$, let
$\Phi^{(q_1,q_2)}$ specify that a character is unramified at $q_1$ and $q_2$,
and  let $\Phi^{(0)}$ make no specification at all,
we find
\begin{align*}
\lim_{X\ra \infty} \frac{F_{D,G}(\Phi^{(q_i)},X)}{F_{D,G}(\Phi^{(0)},X)}&= \frac{1}{\sum_{\chi : \Z_{q_i}^\times \ra G}
 \frac{1}{D(\chi)^{s}}},\\
\lim_{X\ra \infty} \frac{F_{D,G}(\Phi^{(q_1,q_2)},X)}{F_{D,G}(\Phi^{(0)},X)}&= \frac{1}{(\sum_{\chi : \Z_{q_1}^\times \ra G}
 \frac{1}{D(\chi)^{s}})(\sum_{\chi : \Z_{q_2}^\times \ra G}
 \frac{1}{D(\chi)^{s}})}, \textrm{ and}\\
\lim_{X\ra \infty} \frac{F_{D|_{p\Z/p^2\Z},p\Z/p^2\Z}(\Phi^{(0)},X)}{F_{D,\Z/p^2\Z}(\Phi^{(0)},X)}&= 
\lim_{s\ra \frac{1}{p(p-1)}} \frac{F_{D|_{p\Z/p^2\Z},p\Z/p^2\Z}(s,\Phi^{(0)})}{F_{D,\Z/p^2\Z}(s,\Phi^{(0)})}\ne 0,1.
\end{align*}
We can define $D$-probabilities of local specifications for random characters $J/\Q^\times \ra \Z/p^2\Z$ as in 
Equation~\eqref{E:probdef}, essentially replacing the set of surjective characters
$J/\Q^\times \ra G$ by the set of all characters $J/\Q^\times \ra G$.  Then the above tells us that 
$q_1$ ramifying and $q_2$ ramifying are $D$-independent events
for random characters to $\Z/p^2\Z$.
We see that $q_1$ ramifying and $q_2$ ramifying are $D$-independent events
for random characters with image in $p\Z/p^2\Z$.
Also, the probability that a random character to $\Z/p^2\Z$ has image in $p\Z/p^2\Z$ is not 0 or 1.
Since we have $q_i\equiv 1 \pmod{p^2}$, 
there are more maps from $\Z_{q_i}^\times$ to $\Z/p^2\Z$ than to $p\Z/p^2\Z$.
Thus the probabilities of $q_i$ ramifying
in a random character to $\Z/p^2\Z$ and a in random character with image in $p\Z/p^2\Z$ are different.
We have the following simple fact from probability theory.
\begin{lemma}
Let $A$ be an event with positive probability not equal to 1.
 If $E_1$ and $E_2$ are independent, independent given $A$,
and for $i=1,2$ we have that the $\Pr(E_i|A)\ne \Pr(E_i)$, then
$E_1$ and $E_2$ are not independent given not-$A$.
\end{lemma}
So we can conclude that the probabilities of $q_1$ and $q_2$ ramifying
in a random surjective character to $\Z/p^2\Z$, or equivalently in a $\Z/p^2\Z$-extension of $\Q$,
are not independent.
\end{proof}

\begin{proposition}\label{P:notCheb}
Let $q=2,3,5,7,11,$ or $13$.
Given that $q$ is unramified, the discriminant probability that $q$ splits completely
in a $\Z/9\Z$-extension is less than $\frac{1}{9}$.
\end{proposition}

\begin{proof}
From Wright \cite[Theorem I.4]{Wright}, we know that $q$ is unramified with non-zero discriminant probability
in a random
$\Z/9\Z$-extension, and thus is makes sense to formulate the proposition.
First, we let $G=\Z/p^2\Z$ for an arbitrary odd prime $p$.
We let $S=\{q\}$ for some prime $q$, and define $\phi$ on $S$ with $\phi_q$ the trivial character.
We will use the isomorphisms  $\Hom(J/\Q^\times, G)\isom \Hom(\prod_\ell \Z_{\ell}^\times,G) \isom \Hom \left( \left( \prod_{\ell\ne q} \Z_{\ell}^\times \times \Q_q^\times \right)/{\langle q\rangle},G\right)$.
As in Section~\ref{SS:Euler}, for $\epsilon\in \mathcal{A}=\langle q \rangle/\langle q^{p^2} \rangle$
we define
$$
F_{D,G}(s, \phi):= \sum_{\substack{\chi : 
J/\Q^\times \ra G\\ \forall v\in S \ \chi_v=\phi_v }}
\frac{1}{D(\chi)^{s}}
\quad \textrm{and}\quad
  F'_{D,G}(s, \epsilon, \phi):= \sum_{\substack{\chi : 
\prod_{\ell\ne q} \Z_{\ell}^\times \times \Q_q^\times \ra G\\ \forall v\in S \ \chi_v=\phi_v }}
\frac{\zeta_{p^2}^{\chi(\epsilon)}}{D(\chi)^{s}},
$$
and have
  $F_{D,G}(s, \phi)=\frac{1}{|\cA|}\sum_{\epsilon \in \cA }  F'_{D,G}(s,\epsilon, \phi)$.
We have the usual Euler product
\begin{align*}
   F'_{D,G}(s, \epsilon, \phi)
&= \prod_{\ell\ne q} \sum_{\chi_\ell:\Z_l^\times \ra G} \frac{\zeta_{p^2}^{\chi_\ell(\epsilon)}}{D(\chi_\ell)^{s}},
\end{align*}
which has no factor at $q$ because $\phi_q$ is the trivial character.
We see that $F'_{D,G}(s, \epsilon, \phi)$
 only differs from $F_{D,G}(s, \epsilon, \phi')$
(for any choice of $\phi'$ on $S'\ni q$) of Section~\ref{S:anal} by a finite number of factors.
We also see that $F'_{D,G}(s, \epsilon, \phi)/F_{D,G}(s, \epsilon, \phi')$ is entire
and non-zero at $\frac{1}{p(p-1)}$, and thus we conclude
from Lemma~\ref{L:cont} that
$F'_{D,G}(s, \epsilon, \phi)$ can be analytically continued to $\Re(s)\geq\frac{1}{p(p-1)}$
except
for a possible pole  of order at most one at $\frac{1}{p(p-1)}$. From Lemma~\ref{L:whopoles} we have that
 $F'_{D,G}(s, \epsilon, \phi)$
has a pole at $\frac{1}{p(p-1)}$
exactly when $\epsilon^{1/p}\in \Q(\zeta_{p})$, i.e.  when $\epsilon\in\langle q^p \rangle$.

For $\ell\ne q$ and $p\nmid i$,
\begin{align*}
& \sum_{\chi_\ell:\Z_\ell^\times \ra G} \frac{\zeta_{p^2}^{\chi_\ell(q^{pi})}}{D(\chi_\ell)^{s}}=\\
&\begin{cases}
 1 , &\ell\not\equiv 1 \pmod{p};\\
1 +(p-1)\ell^{-(p^2-p)s} , &\ell \equiv 1 \pmod{p}\text{ and }\ell\not\equiv 1 \pmod{p^2};\\
1 +(p-1)\ell^{-(p^2-p)s} -p\ell^{-(p^2-1)s} , &\ell \equiv 1 \text{ or } p\pmod{p^2}\text{ and $q$ not a $p$th power in $\Q_\ell$};\\
1 +(p-1)\ell^{-(p^2-p)s} +(p^2-p)\ell^{-(p^2-1)s} , &\ell \equiv 1 \text{ or } p\pmod{p^2}\text{ and $q$ a $p$th power in $\Q_\ell$}.\\
\end{cases}
\end{align*}
Also,
\begin{align*}
 \sum_{\chi_\ell:\Z_\ell^\times \ra G} &\frac{1}{D(\chi_\ell)^{s}}=\\
&\begin{cases}
 1 , &\ell\not\equiv 1 \pmod{p};\\
1 +(p-1)\ell^{-(p^2-p)s} , &\ell \equiv 1 \pmod{p}\text{ and }\ell\not\equiv 1 \pmod{p^2};\\
1 +(p-1)\ell^{-(p^2-p)s} +(p^2-p)\ell^{-(p^2-1)s} , &\ell \equiv 1 \text{ or } p\pmod{p^2}.\\
\end{cases}
\end{align*}
To find the discriminant probability that a random character to $\Z/p^2\Z$ splits
completely at $q$, given that it is unramified at $q$, we compare 
$F_{D,G}(s, \phi)$ to $F_{D,G}(s, \Phi^{(q)})$ (from the proof of Proposition~\ref{P:notind}),
which counts all characters to $\Z/p^2\Z$, unramified at $q$. 
We have 
\begin{align*}
 F_{D,G}(s, \Phi^{(q)})
=\prod_{\ell\ne q}\left(\sum_{\chi_\ell : \Z_\ell^\times \ra G} \frac{1}{D(\chi_\ell)^{s}} \right) .
\end{align*}
Both $F_{D,G}(s, \phi)$ and $F_{D,G}(s, \Phi^{(q)})$ 
can be meromorphically continued to $\Re(s)\geq\frac{1}{p(p-1)}$, analytic away from $\frac{1}{p(p-1)}$ 
and with a pole of order 1 at
at $\frac{1}{p(p-1)}$.  Thus we can use a Tauberian theorem, as in the end of Section~\ref{S:anal},
 to find that the discriminant probability  of
a random
 character
$\chi :J/\Q^\times \ra \Z/p^2\Z$
being trivial at $q$, given that it is unramified,  is
\begin{align*}
s&=\lim_{s\ra\frac{1}{p(p-1)}}\frac{\frac{1}{|\cA|}\sum_{\epsilon \in \cA }  F'_{D,G}(s,\epsilon, \phi)}{F_{D,G}(s, \Phi^{(q)})}\\&=\frac{1}{p^2}\left( 1 + (p-1)\prod_{\substack{\ell\equiv 1\textrm{ or }p \pmod{p^2}\\ \text{$q$ not a $p$th power in $\Q_\ell$}\\ \ell\ne q}} \frac{(1 + (p-1)\ell^{-1}  -p \ell^{-(p+1)/p})}
{(1 + (p-1)\ell^{-1}  +(p^2-p) \ell^{-(p+1)/p})}  \right).
\end{align*}

\begin{remark}
 Note that $s>\frac{1}{p^2}$ because we know both $F'_{D,G}(s,q^{pi}, \phi)$ and $F_{D,G}(s, \Phi^{(q)})$
do have a pole at $\frac{1}{p(p-1)}$.  Thus we cannot ``resolve'' this proposition by simply considering
all characters $\chi :J/\Q^\times \ra \Z/p^2\Z$ instead of just $\Z/p^2\Z$-extensions.
\end{remark}
We have shown that the discriminant probability of $q$ splitting completely in a random character
$\chi :J/\Q^\times \ra p\Z/p^2\Z$, given that it is unramified at $q$, is $\frac{1}{p}$,
because $D|_{p\Z/p^2\Z}$ is fair and so we can use Corollary~\ref{C:QCheb}.
By the method in the proof of Proposition~\ref{P:notind}, we can compute that
the discriminant probability that a random character
$\chi :J/\Q^\times \ra \Z/p^2\Z$, unramified at $q$, has image in $p\Z/p^2\Z$
is
$$
r=\prod_{\ell\ne q}
\frac{\sum_{\chi : \Z_{\ell}^\times \ra p\Z/p^2\Z}
 \frac{1}{D(\chi)^{1/(p^2-p)}}}{\sum_{\chi : \Z_{\ell}^\times \ra \Z/p^2\Z}
 \frac{1}{D(\chi)^{1/(p^2-p)}}}
=\prod_{\substack{\ell\equiv 1\textrm{ or }p \pmod{p^2}\\ \ell\ne q}} \frac{(1 + (p-1)\ell^{-1})}
{(1 + (p-1)\ell^{-1}  +(p^2-p) \ell^{-(p+1)/p})}.
$$
Thus if $s_1$ is the probability that a random surjective character to $\Z/p^2\Z$ 
is trivial at $q$, given that it is unramified at $q$, we have
$$
s_1=\frac{s-\frac{r}{p}}{1-r}
$$
and thus $s_1>\frac{1}{p^2}$ if and only if $\frac{p^2s-1}{(p-1)r}>1$.
In other words, $s_1>1$ if and only if
$$
\prod_{\substack{\ell\equiv 1\textrm{ or }p \pmod{p^2}\\ q\not\in\Q_\ell^p\\ \ell\ne q}} \frac{(1 + (p-1)\ell^{-1}  -p \ell^{-(p+1)/p})}
{(1 + (p-1)\ell^{-1} )} 
\prod_{\substack{\ell\equiv 1\textrm{ or }p \pmod{p^2}\\q\in\Q_\ell^p\\ \ell\ne q}} \frac{(1 + (p-1)\ell^{-1}  +(p^2-p) \ell^{-(p+1)/p})}
{(1 + (p-1)\ell^{-1} )}
>1 .
$$
We can compute truncations of the above product in PARI/GP \cite{PARI2} for $p=3$, $q=2,3,5,7,11,13$, and $\ell\leq N$, where $N=10^5$
(except when $q=3$ we use $N=10^8$).
We can estimate that the remainder, the product of the terms with $l>N$
is at most
\begin{align*}
 \prod_{\ell >N } (1 + (p^2-p) \ell^{-(p+1)/p)})&\leq \prod_{\ell >N } (1 +  \ell^{-(p+1)/p})^{p^2-p}\\
&\leq \left(1 + \sum_{n >N } n^{-(p+1)/p} \right )^{p^2-p}\\
&\leq \left(1 + pN^{-1/p} \right )^{p^2-p},
\end{align*}
where the sum is over integers $n$.
We can then prove that $s_1 \leq .97$ in all of these cases.
In conclusion, the probability that a random $\Z/9\Z$-extension of $\Q$
splits completely at $q$, given
 that it is unramified at $q$, is less than $\frac{1}{9}$
for $q=2,3,5,7,11,$ or $13$.
\end{proof}

\section{Fair Artin Conductors}\label{S:AC}

For any faithful finite dimensional complex representation $R$ of $G$
and $G$-extension $L$, we have the Artin conductor $C^R(L)$,
which is a counting function (as defined in the beginning of Section~\ref{S:anal}).  
If $R$ is not faithful, then the Artin conductor is not a counting function because it will
have $c^R_G(g)=0$ for non-trivial $g$.
We have seen that for
fair counting functions, the probabilities of local behaviors are nice, but in Section~\ref{S:disc},
we saw that for an example of an unfair counting function, the probabilities are not so well-behaved.
In this section, we give two simple examples of Artin conductors which give fair counting functions.

For a general definition of Artin conductors, see \cite[VII.11]{Neu}.  
The discriminant is given by the Artin conductor of the regular representation.
Since we are only concerned with $G$
abelian, any representation $R$ breaks up as a sum of one dimensional representations,
each of which is determined by the kernel of the action of $G$ on that one dimensional representation.
Suppose $R$ is given by kernels $H_1$,\dots,$H_s$.  Then for $g\in G$, we have $c^R_G(g)=s-\#\{i|g\in H_i\}$.
This can serve as a definition of the Artin conductor at all tame places, which is all 
that concerns fairness.  In other words, for a character $\chi : K_v^\times \ra G$
for $v\nmid |G|$, we have $c^R(\chi)=c^R_G(\chi(y_v))$, where $y_v$ is a generator of tame inertia.
Recall that $m_R$ is the minimum value, other than 0, taken by $c^R_G$, and $\mathfrak{M}=\mathfrak{M}_R={(c^R_G)}^{-1}(m)$.
The counting function is fair if 
$\mathfrak{M} \cap \{g\in G\mid g^r=1 \}$ generates the subgroup $\{g\in G\mid g^r=1 \}$ for all $r$.

We write $G=\prod_i \Z/n_i\Z$, and let $f_i : G \ra \Z/n_i\Z \hookrightarrow \C^*$
be the projection of $G$ to a factor composed with an injection to $\C^*$.
Then $\oplus_i f_i$ gives a fair Artin conductor.  Since the representation is faithful, 
 the Artin conductor of $\oplus_i f_i$ is a counting function.
Also, the elements of $\mathfrak{M}$ are exactly the elements of $G$ that are in all but one $\ker f_i$, and these are the elements with non-zero
coordinates in exactly one factor of $G$.  These elements of $\mathfrak{M}$ generate $G$ in every exponent, and thus the Artin conductor is fair.

Also, $\bigotimes_i f_i \oplus \bigoplus_i f_i$ has a fair Artin conductor.  
We have $\bigcap_i \ker f_i =\{1\}$ and $\ker \left( \bigotimes_i f_i\right) \cap \bigcap_{i\ne j}\ker f_i =\{1\}$,
and so the elements of $\mathfrak{M}$ are exactly the elements of $G$ that are in all but two of the $\ker f_i$
and $\ker \left( \bigotimes_i f_i\right)$.  The elements of $G$ with non-zero
coordinates in exactly one factor are in $\mathfrak{M}$, and they generate $G$ in every exponent, and thus in this case
the Artin conductor is fair.
We can apply these two examples of fair Artin conductors to other factorizations of $G$ into cyclic groups to obtain more examples
of fair Artin conductors.

\section{Further Questions}\label{S:further}
One may ask whether counting abelian extensions by conductor or discriminant is more natural.  
In this paper, we have seen that the probabilities of local behaviors are very nice when counting by conductor and not
so well behaved when counting by discriminant.
While in both cases we can obtain asymptotic counting results for the total number of extensions
(see Section~\ref{S:constant} and \cite{Wright}), in the case of conductor we can express
the constant in the asymptotic count as an Euler product (see Theorem~\ref{T:constant}).  No
Euler product is known for the constant counting abelian extensions by discriminant for a general group $G$
and base field $K$.  So it seems for abelian groups $G$, counting by conductor gives more natural answers.

The other main examples where this global asymptotic counting and computation of local proabilities can be done
are for degree $n$ extensions with Galois closure with group $S_n$ for $n=3,4,5$ (see \cite{S3}, \cite{S4}, \cite{S5}).  In these cases the counting is done by discriminant, and in fact
it is not clear what we might mean by \emph{conductor} in these cases.  Perhaps one should define the conductor to be the 
greatest common divisor of all Artin conductors.  In \cite{S36} the present author and Bhargava count 
these $S_3$ extensions another way; equivalently, we count Galois degree 6 extensions
with Galois group $S_3$ by their discriminant.  In this case, we obtain an asymptotic for the overall count with an Euler product constant and nice local behaviors (simple ratios of probabilities at a given place, and independence at any finite set of places).
In \cite{S36} it is remarked that one can obtain all these nice behaviors for a range of counting functions.

For quartic extensions of $\Q$ having Galois closure with Galois group $D_4$ the overall asymptotic counting by discriminant has been completed (see \cite{D4}),
but
the constant has not been found to have a simple form, and no results for local probabilities analogous to those in this paper have been found.  We wonder if counting these $D_4$ extensions another way would yield nicer results.  In particular, see
\cite[Section 5]{Mass} for a specific counting function one might investigate.

Ellenberg and Venkatesh \cite[Section 4.2]{Ellenberg} suggest that we can try to count extensions of global fields by 
general counting functions (our terminology).
The larger question that is motivated by this paper is which of these counting functions are better than others.
For which counting functions can we obtain an asymptotic total count?  For which counting functions
is the constant in the asymptotic total count an Euler product?  And for 
which counting functions are the local probabilities simple and independent at finite sets of places?
These questions are exactly in line with the questions of Bhargava in \cite[Section 8.2]{B:MF}, except he asks these questions
mainly for counting by discriminant and here we emphasize that the answers will depend on the choice of counting function.


\begin{thebibliography}{99} 

\bibitem{AT}
E. Artin\ and\ J. Tate, {\it Class field theory}, W. A. Benjamin, Inc., New York, 1968. 

\bibitem{S36} 
M. Bhargava\ and\ M. M. Wood, The density of discriminants of $S\sb 3$-sextic number fields, Proc. Amer. Math. Soc. {\bf 136} (2008), no.~5, 1581--1587.

\bibitem{S4} M. Bhargava, The density of discriminants of quartic rings and fields, Ann. of Math. (2) {\bf 162} (2005), no.~2, 1031--1063. 

\bibitem{S5} M. Bhargava, The density of discriminants of quintic rings and fields,  Ann. of Math., to appear.

\bibitem{B:MF} M. Bhargava, Mass formulae for extensions of local fields, and conjectures on the density of number field discriminants, Int. Math. Res. Not. IMRN {\bf 2007}, no.~17, Art. ID rnm052, 20 pp.

\bibitem{Cohen1} H. Cohen, Constructing and counting number fields, in {\it Proceedings of the International Congress of Mathematicians, Vol. II (Beijing, 2002)}, 129--138, Higher Ed. Press, Beijing.


\bibitem{Cohensur}
H. Cohen, F. Diaz y Diaz\ and\ M. Olivier, Counting discriminants of number fields, J. Th\'eor. Nombres Bordeaux {\bf 18} (2006), no.~3, 573--593.

\bibitem{Cohen2} H. Cohen, F. Diaz y Diaz\ and\ M. Olivier, A survey of discriminant counting, in {\it Algorithmic number theory (Sydney, 2002)}, 80--94, Lecture Notes in Comput. Sci., 2369, Springer, Berlin.

\bibitem{D4}
H. Cohen, F. Diaz y Diaz\ and\ M. Olivier, Enumerating quartic dihedral extensions of $\Q$, Compositio Math. {\bf 133} (2002), no.~1, 65--93.

\bibitem{DW} B. Datskovsky\ and\ D. J. Wright, The adelic zeta function associated to the space of binary cubic forms. II. Local theory, J. Reine Angew. Math. {\bf 367} (1986), 27--75.

\bibitem{S3}
H. Davenport\ and\ H. Heilbronn, On the density of discriminants of cubic fields. II, Proc. Roy. Soc. London Ser. A {\bf 322} (1971), no.~1551, 405--420.

\bibitem{DC1}
I. Del Corso\ and\ R. Dvornicich, A converse of Artin's density theorem: the case of cubic fields, J. Number Theory {\bf 45} (1993), no.~1, 28--44.

\bibitem{DC2}
I. Del Corso\ and\ R. Dvornicich, Uniformity over primes of unramified splittings, Mathematika {\bf 45} (1998), no.~1, 177--189.

\bibitem{DC3}
I. Del Corso\ and\ R. Dvornicich, Uniformity over primes of tamely ramified splittings, Manuscripta Math. {\bf 101} (2000), no.~2, 239--266.

\bibitem{Ellenberg}  J. S. Ellenberg\ and\ A. Venkatesh, Counting extensions of function fields with bounded discriminant and specified Galois group, in {\it Geometric methods in algebra and number theory}, 151--168, Progr. Math., 235, Birkh\"auser, Boston, Boston, MA.

\bibitem{Grunwald}
W. Grunwald, Ein allgemeines Existenztheorem f\"{u}r algebraische Zahlk\"{o}rper,
J. Reine Angew. Math. {\bf 169} (1933) 103--107.

\bibitem{Malle}
G. Malle, On the distribution of Galois groups. II, Experiment. Math. {\bf 13} (2004), no.~2, 129--135.

\bibitem{MakiDisc}
S. M\"aki, On the density of abelian number fields, Ann. Acad. Sci. Fenn. Ser. A I Math. Dissertationes No. 54 (1985), 104 pp.

\bibitem{MakiCond}
S. M\"aki, The conductor density of abelian number fields, J. London Math. Soc. (2) {\bf 47} (1993), no.~1, 18--30.

\bibitem{N}
Narkiewicz, W., translated by S. Kanemitsu,
\textit{Number Theory,} 
 World Scientific,  1983.

\bibitem{Neu}
J. Neukirch, {\it Algebraic number theory}, Translated from the 1992 German original and with a note by Norbert Schappacher, Springer, Berlin, 1999.

\bibitem{PARI2}
PARI/GP, version {\tt 2.3.2}, Bordeaux, 2006, {\tt http://pari.math.u-bordeaux.fr/}.

\bibitem{Taylor} M. J. Taylor, On the equidistribution of Frobenius in cyclic extensions of a number field, J. London Math. Soc. (2) {\bf 29} (1984), no.~2, 211--223.

\bibitem{vdP}
C. E. van der Ploeg, On a converse to the Tschebotarev density theorem, J. Austral. Math. Soc. Ser. A {\bf 44} (1988), no.~3, 287--293.

\bibitem{Wang}
S. Wang, On Grunwald's theorem, Ann. of Math. (2) {\bf 51} (1950), 471--484.

\bibitem{Mass} M. M. Wood, Mass formulas for local Galois representations to wreath products and cross products, Algebra Number Theory {\bf 2} (2008), no.~4, 391--405.

\bibitem{GWex} M. M. Wood, $P$-adic algebras that do not occur as completions of number fields, in preparation.

\bibitem{Wright}
D. J. Wright, Distribution of discriminants of abelian extensions, Proc. London Math. Soc. (3) {\bf 58} (1989), no.~1, 17--50. 

\end{thebibliography}
\end{document}